\documentclass[12pt]{amsart}

\usepackage{amssymb,amsmath,amscd,amsfonts,a4,psfrag,graphicx,
latexsym,epsfig,color}
\usepackage[active]{srcltx}
\usepackage[all]{xy}
\usepackage{stmaryrd}
\everymath{\displaystyle}

\newtheorem{teo}{Theorem}[section]

\newtheorem{defi}[teo]{Definition}

\newtheorem{proposition}{Proposition}
\newtheorem{theorem}{Theorem}

\newtheorem{lemma}{Lemma}

\title[Characters
]{Evaluation of characters of smooth representations of $GL(2,{\mathcal O}):$\\
 I.Strongly primitive representations of even level.}

\author{Ph.Roche}

\address{ IMAG, Univ Montpellier, CNRS, Montpellier, France}
\email{philippe.roche@univ-montp2.fr}

\begin{document}

\today

\maketitle

\begin{abstract}  
Let $F$ be a local field,  let ${\mathcal O}$ be its integer ring and $\varpi$ a uniformizer of its maximal ideal.  
To an irreducible complex finite dimensional smooth representation $\pi$ of $GL(2,{\mathcal O})$ is associated a pair of positive integers $k, k'$ called the level and the sublevel of $\pi.$
The level  is the smallest integer $k$ such that $\pi$ factorizes through  the finite group 
$GL(2,{\mathcal O}/\varpi^k {\mathcal O})$, whereas the sublevel  is the smallest integer $k'\leq k$ such that there exists $\chi,$ one dimensional representation of  $GL(2,{\mathcal O}),$ such 
that $\pi\otimes \chi$ factorizes through the finite group  $GL(2,{\mathcal O}/\varpi^{k'} {\mathcal O}).$
A representation  of $GL(2,{\mathcal O})$ is said strongly primitive  if the level and  sublevel are equal.
The classification of  smooth finite dimensional representations  of $GL(2,{\mathcal O})$ is equivalent to the classification of strongly primitive irreducible representations  of $GL(2,{\mathcal O}).$

In this first article we describe explicitely the even level strongly primitive irreducible finite dimensional  complex representations of 
$GL(2,{\mathcal O})$ along the lines of \cite{St} and  \cite{LCW} using Clifford theory. 
In the case where the characteristic $p$ of the residue field is not equal to $2,$ we give exact formulas
 for the characters of these representations in most cases by reducing them to the evaluation of Gauss sums, Kloosterman sums and Sali\'e sums for the finite ring ${\mathcal O}/\varpi^k{\mathcal O}.$ It generalizes the   work of \cite{LCW} which was devoted to $F={\mathbb Q}_p.$
The second article  \cite{BaRo} will give  the  evaluation  of characters in the  odd level case and the exact expressions  for certain  generalized Zeta function representations  \cite{Ro}  of 
$PGL(2,{\mathcal O})$.

\end{abstract}
\medskip


\tableofcontents

\section{Introduction}
Let $F$ be a local field, ${\mathcal O}$ its ring of integers and $p$ the characteristic of the residue field. $GL(2,{\mathcal O})$ is a maximal compact subgroup of $GL(2,F)$ and it is a fundamental theorem that the smooth irreducible complex representations of $GL(2,F)$ are admissible i.e decompose with finite multiplicities in terms of smooth irreducible representations of $GL(2,{\mathcal O}).$  It is therefore of interest to have a classification of irreducible smooth complex representations of  $GL(2,{\mathcal O})$   and to have closed expressions for the characters of this group.

Classification of irreducible smooth complex representations of  $GL(2,{\mathcal O})$ has been obtained in \cite{St}  using Clifford theory \cite{Is}.
Finer results such 
as the explicit evaluation of characters has been obtained in \cite{LCW} only when $F={\mathbb Q}_p$ and $p\not=2$. 

In this article we generalize this last work for arbitrary local field $F$ with $p\not=2$ using similar methods. 

In section 2, we recall the classification of conjugacy classes of  $GL(2,{\mathcal O})$ following \cite{AOPV}. When $p\not=2$ we give a classification of them in  the form given in \cite{LCW}.

In section 3, we recall   the classification of smooth irreducible representations  of $GL(2,{\mathcal O})$ which are strongly primitive of even level  following  \cite{St}  and give a detailed and simpler  description of these representations. 

In section 4, 
 we evaluate the characters of these representations using Frobenius formula.  The expression of the  characters can in most cases be evaluated in  a closed form by reducing them to twisted Kloosterman sums associated to the finite ring ${\mathcal O}/\varpi^k{\mathcal O}.$ Note that, not to diminish the value of the work  \cite{LCW}, we have simplified and sometimes 
corrected their work.

The original motivation for our   work was the evaluation of generalized Zeta representation function \cite{Voll1, JZ, Ro} of the group $PGL(2,{\mathcal O}).$
In order to keep the length of the present  article reasonable, we have computed the characters of the representations only for those which are strongly primitive of even level.
In a forthcoming article \cite{BaRo} we  evaluate the characters of representations which are strongly primitive of odd length and by mixing these two results we  give closed expression of the evaluation  of  certain generalized Zeta representation functions of $PGL(2,{\mathcal O})$

\section{Conjugacy classes}
Let ${\mathcal  A}$ be  a local principal ring, let ${\mathfrak M}$ be the maximal ideal of
 ${\mathcal A},$ $\varpi$ a uniformiser of $ {\mathfrak M},$ $\Bbbk$ the  residual field of characteristic $p.$ Let $r\in {\mathbb N}\cup \{\infty\}$ be the length of 
${\mathcal A},$ i.e the smallest positive integer  $r,$  if it exists, such that ${\mathfrak M}^r=\{0\}$,  if not  we define $r=\infty.$ By convention we denote $\varpi^\infty=0,$ and $ {\mathfrak M}^\infty=\{0\}.$
If $r\not=\infty, \llbracket 0,r \rrbracket=\{0,\cdots,r\}$ and $\llbracket 0, \infty \rrbracket={\mathbb N}\cup\{\infty\}.$
We recall the classification of similarity classes of matrices of $M_2({\mathcal A})$ as given in \cite{AOPV} (Theorem 2.2). For $i\in\llbracket 0,r \rrbracket$  we denote ${\mathcal A}_i={\mathcal A}/{\mathfrak  M}^i, $ and for $i\not=\infty$ we choose $s_i: {\mathcal A}_i\rightarrow {\mathcal A }$  sections of  the canonical projections $p_i:{\mathcal A}\rightarrow {\mathcal A}_i.$  
 We denote ${\mathbb A}_i\subset {\mathcal A}$ the image of ${\mathcal A}_i$ under $s_i$ for $i\not=\infty.$ We choose $s_i$ such that ${\mathbb A}_0=\{0\}$,   $s_1(0)=0,$ and for  $0<i<\infty$, ${\mathbb A}_i=\{\sum_{j=0}^{i-1} a_j \varpi^j, a_j\in {\mathbb A}_1 \}$, ${\mathbb A}_i$ is in bijection with the set ${\mathcal A}_i.$  Let $j\in  \llbracket 0,r \rrbracket$, $j\not=\infty,$ we denote $\rho_j: {\mathcal A}\rightarrow  \varpi^j {\mathcal A}, a\mapsto\varpi^j a,$ which after quotienting by the kernel defines an isomorphism of ${\mathcal A}$-module  $ {\bar\rho}_j:{\mathcal A}_{r-j} \xrightarrow {\sim } \varpi^j{\mathcal A}.$ 

The following easy lemma is central in the classification of \cite{AOPV}. 

\begin{lemma}
Let  $X\in M_2(\mathcal A)$, it can be written as $X=A+\varpi^j B$ where $j\in\llbracket 0,r \rrbracket $ is maximal such that $X$ is congruent modulo  ${\mathfrak M}^j$ to a scalar matrix $A=\alpha I,\alpha\in {\mathcal A}.$  If $j=\infty$ then $X=\alpha I, \alpha\in {\mathcal A}.$
If  $j\not=\infty$, $\alpha$ can be chosen in  ${\mathbb A}_j$ and is unique,  $B$  is unique mod ${\mathfrak  M}^{r-j}$  and moreover is a cyclic matrix i.e there exist $v\in {\mathcal A}_{r-j}$ such that $(v, Bv)$ is a basis of ${\mathcal A}_{r-j}^{2}.$ 
\end{lemma}
Note that $\alpha=0$ when $j=0.$ 

If $(a,b)\in {\mathcal A}$ we denote 
$C(a,b)=\begin{pmatrix}0&1\\a&b\end{pmatrix}.$

As a result one obtains the theorem (theorem 2.2 of \cite{AOPV}):

\begin{theorem}\label{theoclassificationconjugyclass}
Let $X\in M_2({\mathcal A})$ and $j, \alpha, B$  associated to $X$ by the previous lemma, then $X$ is similar to the matrix
 $\alpha I+\varpi^j C(-det(B),tr(B))=\begin{pmatrix}\alpha&\varpi^j\\-\varpi^j det(B)&\alpha+\varpi^j tr(B)
\end{pmatrix}.$ 
Inversely,  given  $j\in\{0,...,r\}$, $j\not=\infty,$ $\alpha\in {\mathbb A}_j$ , and a couple $(\varpi^j\beta,\varpi^j\gamma)\in  {\bar\rho}_j({\mathcal A}_{r-j}),$  there exists a unique class of similarity matrix having $\begin{pmatrix}\alpha&\varpi^j\\-\varpi^j \gamma&\alpha+\varpi^j \beta\end{pmatrix}$ as representative. If $j=\infty$ and $\alpha\in  {\mathcal A},$ the class of similarity matrix having $\alpha I$ as representative consists only on this matrix.

\end{theorem}
As a result the conjugacy classes of $GL(2,{\mathcal  A})$ are in bijection with the subset of these representatives defined by the additional condition that the determinant is invertible. This last condition can also be written: if $j=0$ then $\alpha=0$ and $\gamma=det(B)\in {\mathcal  A}^\times $  and if $j\geq 1$ then $\alpha\in {\mathcal  A}^\times $.

A further  classification, simpler,  is obtained when  the characteristic of the residual field $\Bbbk$ is different of $2.$ 
\begin{proposition}\label{propclassificationconjugyclass}
A set of representatives of similarity classes of $M_2({\mathcal A})$ are given by the set of matrices 
 $\begin{pmatrix}\alpha&\varpi^j\\\varpi^j \beta &\alpha\end{pmatrix}, $  $j\in \llbracket 0,r \rrbracket,$ $j\not=\infty,$  $\alpha\in {\mathcal A}$, $\varpi^j\beta\in {\bar\rho}_j({\mathcal A}_{r-j})$ with the addition of the case $j=\infty$ (when $r=\infty)$ and the set of matrices $\alpha I$ with $\alpha\in {\mathcal A}.$
\end{proposition}

\proof 
Let $y\in {\mathcal A}, $ we denote $Y=\begin{pmatrix}1&0\\y&1
\end{pmatrix},$  we have $Y(\alpha I+\varpi^jC(-\gamma,\beta) )Y^{-1}=
\begin{pmatrix}\alpha-\varpi^j y&\varpi^j\\ \varpi^j(-\gamma-y^2-\beta y))& \alpha+\varpi^j\beta+\varpi^j y\end{pmatrix}.$
Therefore if $2$ is invertible in ${\mathcal A},$ we can choose $y=-\frac{1}{2}\beta$ in order to impose that  the elements on the diagonal are equal.
\endproof

Remark: The proposition \ref{propclassificationconjugyclass} is easily shown to be false when $p=2$. Indeed fix $j=0$ and $\beta$ invertible. The matrix $(\alpha I+\varpi^jC(-\gamma,\beta) )$ has its trace equal to $2\alpha+\beta$, therefore it cannot be similar to a matrix having  the same elements on the diagonal which  trace, multiple of $2,$ is therefore non invertible. 

Let $F$ be local field, we assume that the characteristic $p$  of the residual field  is different of $2$  and let ${\mathcal O}$ be its integer ring. We will now use the previous classification when ${\mathcal  A}={\mathcal O}$ and  the length $r=\infty.$ One obtains a generalisation of the classification obtained in \cite{LCW} for the case ${\mathcal A}=
{\mathbb Z}_p.$
Let fix $\epsilon\in {\mathcal O}^{\times} $ which is not a square, it always exists because  $p\not=2.$ 
\begin{proposition} \label{theoclassificationconjugyclasstwo}

A set of representatives of conjugacy classes of $GL(2,{\mathcal O})$ is given by:

\begin{itemize}
\item 
$I_\alpha=\begin{pmatrix}\alpha&0\\ 0&\alpha \end{pmatrix}, \alpha\in{\mathcal O}^\times$
(Scalar class)\item

$B_{i,\alpha,\beta}=\begin{pmatrix}\alpha&\varpi^{i+1}\beta\\ \varpi^{i}&\alpha\end{pmatrix}, i\in {\mathbb N}, \alpha\in{\mathcal O}^\times, \beta\in {\mathcal O}$ (Unipotent class)
\item
$C_{i,\alpha,\beta}=\begin{pmatrix}\alpha&\varpi^i \epsilon \beta\\ \varpi^i \beta&\alpha \end{pmatrix}, i\in {\mathbb N}, \alpha\in{\mathcal O}, \beta\in {\mathcal O}^\times, \alpha^2-\epsilon \beta^2\varpi^{2i}\in {\mathcal O}^\times$
(Elliptic class)
\item

$D_{i,\alpha,\delta}=\begin{pmatrix}\alpha&0\\ 0&\delta \end{pmatrix}, i\in {\mathbb N}, \alpha\in {\mathcal O}^\times, \delta\in{\mathcal O}^\times,\alpha-\delta\in \varpi^i{\mathcal O}^{\times}$
(Diagonal class).
\end{itemize}
\end{proposition}

\proof
One uses the proposition (\ref{propclassificationconjugyclass}) giving a set of representatives of conjugacy classes of  $GL(2,{\mathcal O})$ to be $\begin{pmatrix}\alpha&\varpi^j\\\varpi^j \beta&\alpha\end{pmatrix}, $  $ j\in {\mathbb N}$, $\alpha\in {\mathcal O}$, $\beta\in {\mathcal O}$ with $\alpha^2-\varpi^{2j}\beta\in {\mathcal O}^\times,$ with the addition of the matrices 
 $\alpha I, $ $\alpha\in {\mathcal O}^{\times}$ corresponding to $j=\infty.$

Let  $j\in {\mathbb N}$, if $\beta$ is not invertible in ${\mathcal O}$ then $\beta=\varpi \beta'$ and a representative of this conjugacy class  is given by $B_{j,\alpha,\beta'}.$
If $\beta$ is invertible, there are two possibilities: it is a  square or not.
If $\beta=\mu^2,$  let  $P=\begin{pmatrix} \mu&-\mu\\ 1&1
\end{pmatrix}$, $P$ is invertible ($det(P)=2\mu$) when $p\not=2$ and we have  $P^{-1}\begin{pmatrix} \alpha &\varpi^j\mu^2\\ \varpi^j& \alpha\end{pmatrix}P=
\begin{pmatrix} \alpha+\varpi^j\mu &0\\ 0& \alpha-\varpi^j\mu\end{pmatrix}.$  Therefore a representative of this conjugacy class is given by the matrix $D_{j,\alpha+\varpi^j\mu,\alpha-\varpi^j\mu}.$
We have fixed $\epsilon\in {\mathcal O}^{\times} $ which is not a square, therefore if $\beta$ is not a square $\beta \epsilon^{-1}$ is a square $\nu^2$. This comes from the fact that by Hensel lemma an invertible element is a square in ${\mathcal O}$ if and only if it is a square mod ${\mathfrak M}.$
If we denote $P=\begin{pmatrix} 1&0\\ 0&\nu
\end{pmatrix}$ we have $P \begin{pmatrix} \alpha &\varpi^j\epsilon \nu^2\\ \varpi^j& \alpha\end{pmatrix} P^{-1}=
\begin{pmatrix} \alpha &\varpi^j\epsilon \nu\\ \varpi^j\nu& \alpha\end{pmatrix}=C_{j,\alpha,\nu}.$
As a result the set of matrices defined in the proposition is a set of representatives of the conjugacy classes of $GL(2,{\mathcal  O}).$

The name of the classes comes from the name of the projection  of the matrix in $GL(2, {\Bbbk}).$

\endproof


\section{Irreducible finite dimensional complex smooth representations of $GL(2,\mathcal O)$}
In this section $F$ is a local field, $v$ the additive valuation normalized by $v(\varpi)=1,$ ${\mathcal O}$ is the ring of integers of $F$ and $p$  the characteristic of the residual field $\Bbbk.$ We denote $q$ the cardinal of $\Bbbk.$  We dot not assume  in this section, unless explicitely stated, that $p\not=2.$

Let $r\in {\mathbb N}_{>0}$, we denote ${\mathcal O}_r={\mathcal O}/\varpi^r {\mathcal O},$ and we define $G^{(r)}=GL(2,{\mathcal O}_r).$ $GL(2, {\mathcal O})$ is the profinite group $\varprojlim G^{(r)}$.  We denote $p_r:GL(2, {\mathcal O})\rightarrow GL(2,{\mathcal O}_r)$ the canonical maps.

\begin{defi} If $\pi$ is a finite dimensional complex  smooth representation of $GL(2, {\mathcal O})$ then there exists an integer $k$ such that $\pi$ factorizes through $p_k$ as $\pi=\pi_k\circ p_k$ where $\pi_k$ is a representation of $GL(2, {\mathcal O}_k)$. $\pi$ is irreducible if and only if $\pi_k$ is.

The smallest of these $k$  is by definition the {\bf level} of $\pi.$

A representation of $G^{(r)}$ of level $r$ is said to be {\bf primitive}.
\end{defi}

As a result the classifications of irreducible finite dimensional smooth representations of $GL(2,{\mathcal O})$ of level less than $r$ is equivalent to the classification of irreducible  finite dimensional complex representations of the finite group $GL(2,{\mathcal O}_r)$.  

Remark 1: Note that the theorem  \ref{theoclassificationconjugyclass} applies as well when  ${\mathcal  A}={\mathcal O/\varpi^r {\mathcal O}}$ where $r$ is any positive integer.
Therefore a set of representative of conjugacy classes of  $GL(2,\mathcal O/\varpi^r {\mathcal O})$  is given by 
$
\{\alpha I+C(-\varpi^j \beta', \varpi^j \alpha'), \alpha\in {\mathbb A}_j, j=0,\alpha'=0,\beta'\in {{\mathcal O}_r}^\times \text{ or } 1\leq j\leq r, \varpi^j \beta', \varpi^j \alpha'\in \bar{\rho}_j({{\mathcal O}}_{r-j})\}.
$
The cardinal of this set is $q^{r-1}(q-1)+\sum_{j=1}^{r} (q-1)q^{j-1}q^{r-j}q^{r-j}=q^{r-1}(q^{r+1}-1).$ Therefore the number $n_r$ of conjugacy classes of $GL(2,\mathcal O/\varpi^r {\mathcal O})$   is equal to $n_r=q^{r-1}(q^{r+1}-1)$ which is also the number $a_r$ of irreducible finite dimensional complex representations up to isomorphism of the finite group $GL(2,\mathcal O/\varpi^r {\mathcal O})$. As a result the number $b_r$ of  irreducible finite dimensional complex smooth representations up to isomorphism of level $r$ of $GL(2,{\mathcal O })$ is $b_r=n_r-n_{r-1}.$

Remark 2: When $p\not=2.$, we will use the following classification of conjugacy classes of  $GL(2,{\mathcal O}_r)$, which is  a direct application of proposition  (\ref{theoclassificationconjugyclasstwo}).
Let fix $\epsilon\in {\mathcal O}_r^{\times} $ which is not a square, it always exists because  $p\not=2.$ 
\begin{proposition} \label{theoclassificationconjugyclassthree}

A set of representatives of conjugacy classes of $GL(2,{\mathcal O_r})$ is given by:

\begin{itemize}
\item 
$I_\alpha=\begin{pmatrix}\alpha&0\\ 0&\alpha \end{pmatrix}, \alpha\in{\mathcal O}_r^\times$
\item

$B_{i,\alpha,\beta}=\begin{pmatrix}\alpha&\varpi^{i+1}\beta\\ \varpi^{i}&\alpha\end{pmatrix}, i\in \llbracket 0,r-1\rrbracket, \alpha\in{\mathcal O}_r^\times, \beta\in {\mathcal O}_r$
\item
$C_{i,\alpha,\beta}=\begin{pmatrix}\alpha&\varpi^i \epsilon \beta\\ \varpi^i \beta&\alpha \end{pmatrix}, i\in \llbracket 0,r-1\rrbracket, \alpha\in{\mathcal O}_r, \beta\in {\mathcal O}_r^\times, \alpha^2-\epsilon \beta^2\varpi^{2i}\in {\mathcal O}_r^\times$

\item

$D_{i,\alpha,\delta}=\begin{pmatrix}\alpha&0\\ 0&\delta \end{pmatrix}, i\in  \llbracket 0,r-1\rrbracket , \alpha\in {\mathcal O}^\times, \delta\in{\mathcal O}^\times,\alpha-\delta\in \varpi^i{\mathcal O}^{\times}.$
\end{itemize}
\end{proposition}

The problem of classifying the irreducible finite dimensional complex representations  of $GL(2, {\mathcal O}_r)$ can be completely understood and in great detail using Clifford theory, this is what we review in the sequel.

For $0\leq i\leq r,$  let $K_i^{(r)}=\{g\in G^{(r)}, g=I\, \text{mod}\, {\varpi^i}\}.$ 
If the context is clear we will forget the upper index $r.$

We have $\{I\}=K_r^{(r)}\subset K_{r-1}^{(r)}\subset ...\subset K_0^{(r)}=G^{(r)}.$
The isomorphism $\bar{\rho}_j:{\mathcal O}_{r-j} \to \varpi^j{\mathcal O}_r$ is extended to an isomorphism 
$\bar{\rho}_j:
M_2({\mathcal O}_{r-j} )\to  \varpi^jM_2({\mathcal O}_r).$ 
Having fixed a set of compatible section $s_j$ of ${\mathcal O}_j$ (like in section 1)  for $0\leq j\leq r$, we denote ${\mathbb O}_j$ the image of $s_j,$ ${\mathbb O}_j$ and ${\mathcal O}_j$ are in bijection. We denote  ${\mathbb O}_j^{\times}$ the invertible elements of  ${\mathbb O}_j.$

The following properties hold:
\begin{proposition}
\begin{enumerate}
\item
$K_{i}^{(r)}$ is a normal subgroup of $G^{(r)}.$
\item
$G^{(r)}/K_{i}^{(r)}$ is isomorphic to $G^{(i)}, i>0.$

\item 
$K_{i}^{(r)}=I+\varpi^i M_2({\mathcal O}_r)$ if $i>0.$

\item $K_{i}^{(r)}$ is abelian if $i\geq r/2$ and  $(M_2({\mathcal  O}_{r-i}),+)\rightarrow K_{i}^{(r)}, x\mapsto I+\bar{\rho}_i(x)$  is an isomorphism of abelian group if $i\geq r/2$ where  $M_2({\mathcal  O}_{r-i})$ is endowed with the addition of matrix group law. 
\item $\vert K_i^{(r)}\vert =q^{4(r-i)}$ if $i>0.$
\item $\vert G^{(r)}\vert =q^{4r-3}(q+1)(q-1)^2, r\geq 1.$
\end{enumerate}
\end{proposition}

\proof  The only nontrivial result is the computation of $\vert G^{(r)}\vert$. Let $p:GL(2,{\mathcal O}_r)\rightarrow GL(2,{\Bbbk})$ the canonical map, the kernel of $p$ is 
$K_1^{(r)}=I+\varpi M({\mathcal O}_r)$, which cardinal is $q^{4(r-1)}.$ We have $\vert GL(2,{\mathbb F}_{q})\vert=q(q+1)(q-1)^2$ from which the proposition follows.
\endproof

We will define $l=\lfloor \frac{r+1}{2}\rfloor $ and $l'=\lfloor \frac{r}{2}\rfloor .$ We have $l+l'=r$ and $l$ is the smallest  integer $i$ with $i\geq r/2.$

Let us fix a smooth  additive character $\psi^{(r)}:({\mathcal O},+)\rightarrow {\mathbb C}^\times$ of level $r$ which means that   the kernel of $\psi^{(r)}$ contains ${\mathfrak M}^r$ but not
 ${\mathfrak M}^{r-1}$ (such character always exists).  If the context is clear we will denote it simply by $\psi.$

We first recall a simple description of  characters of the abelian groups $K_i^{(r)}$ for $i\geq r/2.$ Let $\beta \in M_2({\mathcal O}_{r})$, one defines $\psi_{\beta}:K_{i}^{(r)}\rightarrow {\mathbb C}^{\times }$ by 
$\psi_{\beta}(x)=\psi^{}(Tr(\beta(x-I))).$ $\psi_\beta$ depends only on  $\varpi^i\beta,$  therefore the map  $M_2({\mathcal O}_{r})\rightarrow Hom(K_{i}^{(r)},{\mathbb C}^\times), \beta\mapsto \psi_{\beta}$ factorizes through an isomorphism (because $\psi^{(r)}$ is of level $r$) $M_2({\mathcal O}_{r-i})\xrightarrow {\sim } Hom(K_{i}^{(r)},{\mathbb C}^\times),\beta\mapsto \psi_{\beta}.$

We use the following theorem of Clifford theory \cite{Is} recalled in \cite{St} (Theorem 2.1):
let $G$ be a finite group and $N$ a normal subgroup of $G$. $G$ acts on the set of representations of $N$ by conjugation: if $\rho$ is a representation of $N$ then for $g\in G$ we denote $\rho^g$ the representation of $N$ defined  by $\rho^g(n)=\rho(gng^{-1}).$ For any irreducible representation $\rho$ of $N$, we define the stabilizer $T(\rho)$ as being the subgroup of $G$ defined by $T(\rho)=\{g\in G, \rho^g \text{ is isomorphic to }\rho\},$ $T(\rho)$ always contains $N.$
Assume that $\rho$ is an irreducible representation of $N,$  then the set of irreducible representations of $G$ which restriction to $N$ contains $\rho$ is in bijection with the set of irreducible representations of $T(\rho)$ which restriction to $N$ contains $\rho.$ More precisely, if $A=\{\theta\in Irr(T(\rho)), Res^{T(\rho)}_N(\theta)\; \text{contains} \;\rho \}$ and $B=\{\pi\in Irr(G), Res^{G}_N(\pi)\; \text{contains} \;\rho \},$ then $\theta\mapsto Ind^G_{T(\rho)}(\theta)$ is a bijection from $A$ to $B.$
Moreover if $\pi$ is an irreducible representation of $G$ then $Res^G_N(\pi)=e(\bigoplus_{\rho\in \Omega}\rho)$ where $\Omega$ is an orbit of the action of $G$ on the set of classes of irreducible representations of $N$ and $e$ is a positive integer.

We will use this theorem and apply it to $G=GL(2,{\mathcal O}_r)$ and $N=K^{(r)}_l$ with $l$ the smallest integer greater than $r/2$, the reason being that $K_l$ is abelian and therefore the irreducible representations of  $K_l$ are one dimensional  and simple to describe and  if $\rho$ is a one dimensional representation of $K_l$ the condition that $l$ is the smallest implies that the stabilizer of $K_l$ is bigger than $K_l$ but not too much.
We will heavily use the method of \cite{St} but we will be more precise in the description of the representations of the stabilisers. This is important for the computation of the characters of 
$GL(2,{\mathcal O}_r).$

$G^{(r)}$ acts on $Hom(K_l,{\mathbb C}^\times)\simeq M_2({\mathcal O}_{l'}) $,   the orbits are  analysed according to their reductions mod ${\mathfrak M}.$ 

In $M_2({\Bbbk})$ there are $4$ types of similarity equivalence classes:
\begin{itemize}
\item  type $c_1$ (scalar): $\begin{pmatrix}
a&0\\0&a
\end{pmatrix}, a\in \Bbbk$
\item type $c_2$ (diagonal):
$\begin{pmatrix}
a&0\\0&d
\end{pmatrix} ,a,d\in \Bbbk,a\not=d$
\item type $c_3$ (elliptic):
$\begin{pmatrix}
0&1\\-\Delta&s
\end{pmatrix}, \Delta,s\in \Bbbk, x^2-sx+\Delta $ irreducible in $\Bbbk[x]$
\item type $c_4$ (unipotent):
$\begin{pmatrix}
a&1\\0&a
\end{pmatrix}, a\in \Bbbk.$
\end{itemize}

If $a\in {\mathcal O}_{l'}$ we denote $\bar{a}\in\Bbbk$ its reduction mod ${\mathfrak M}$, if $\beta \in M_2({\mathcal O}_{l'})$ we denote $\bar{\beta}\in M_2(\Bbbk)$  its reduction mod ${\mathfrak M}.$ 
If $\beta\in M_2({\mathcal O}_{l'}), \psi_{\beta}\vert_{K_{r-1}^{(r)}}$ depends only on $\bar{\beta}$, we denote it $\psi_{\bar{\beta}}.$

Let $\pi$ be an irreducible representation of $G^{(r)}$ acting on $V$,  $\pi\vert_{K_l}$ decomposes as a direct sum of one dimensional representations and we have

$\pi\vert_{K_l^{(r)}}=e\bigoplus_{\beta\in \Omega} \psi_{\beta}$ where  $\Omega$ is an orbit under $G^{(r)}.$ 
We have $\pi\vert_{K_{r-1}^{(r)}}=e\bigoplus_{\beta\in \Omega} \psi_{\bar{\beta}}.$
All these ${\bar{\beta}} $ are in the same orbit under the action of $GL(2,\Bbbk).$
Therefore we can distinguish two cases:
all $\bar{\beta}$ are  the nul  matrix or none of them are zero.
 In the first case  this means that $\pi\vert_{K_{r-1}^{(r)}}$ is a direct sum of trivial representations, therefore $\pi$ factorises as $\pi:G^{(r)}\rightarrow G^{(r)}/{K_{r-1}^{(r)}}\rightarrow GL(V)$, which means, after using $G^{(r)}/{K_{r-1}^{(r)}}\simeq G^{(r-1)}$ that $\pi$ is of level less or equal to  $r-1.$
If, on the contrary , one (all) of the  $\bar{\beta}$ is not the nul  matrix then $\pi$ is of level $r.$

The case $\bar{\beta}$ is of type $c_1$ i.e $\bar{\beta}=a I$ with $a\not=0$ is interesting.
Let $x\in K_{r-1}^{(r)}$ we have
 $\psi_{a I}(x)=\psi^{(r)}(a(Tr(x-I)))=\psi^{(r)}(a(det(x)-1))=\chi_{a}\circ det(x)$ where $\chi_a$ is a character $K_{r-1}^{(r)}\rightarrow {\mathbb C}^{\times}.$ From the theory of extension of characters of abelian group, $\chi_a$ can be extended to a character $\tilde{\chi}_a:{\mathcal O_{r}}^\times\rightarrow {\mathbb C}^\times.$ 
We will denote $\tilde{\psi}_a$  a one dimensional representation of $G^{(r)}$ extending $\psi_{aI}$   by $\tilde{\psi}_a={\tilde{\chi}_a}\circ det.$
Note that $\tilde{\psi}_a$ is of level $r$ because $a\not=0.$
The representation $\pi$ satisfies $\pi=\tilde{\psi}_a\otimes\pi'$ where $\pi'$ is of level less or equal to $r-1.$
This result  motivates the introduction of the notion of sublevel of a complex finite dimensional irreducible smooth representation $\pi$ of $GL(2, {\mathcal O}).$ 
\begin{defi}
The sublevel  is the smallest integer $k$ such that there exists $\chi,$ one dimensional representation of  $GL(2,{\mathcal O}),$ such 
that $\pi\otimes \chi$ factorizes through the finite group  $GL(2,{\mathcal O}_{k}).$
Because the level always exists,  the sublevel always exists and is  less or equal to the level.
A representation of $GL(2, {\mathcal O})$ of level k which sublevel is also $k$ will be called {\bf strongly primitive of level $k$}.
\end{defi}
 Let $n_r$ the the number of conjugacy classes of $G^{(r)},$ let $a_r$ the number of non isomorphic irreducible representations of $G^{(r)}$,  $b_r $ $(resp \;b'_r) $ the number of non isomorphic primitive (strongly primitive)  representations of $G^{(r)}$. From the discussion above  we have $b_r=b_r'+(q-1)a_{r-1}$ which implies $n_r-qn_{r-1}=b_r'.$

Finally one obtains the following proposition \cite{St}:
\begin{proposition}\label{classificationorbitbeta}
Let $\pi$ be an irreducible representation of $G^{(r)}$ and let $\beta$ be an  element in the orbit $\Omega$ of the decomposition of $\pi\vert_{K_l^{(r)}},$
then  $\beta$ is conjugated  under $G^{(r)}$ to one of these elements:
\begin{itemize}
\item ($C_1$) $\begin{pmatrix}a&0\\0&d
\end{pmatrix}, a,d\in {\mathbb O} _{l'}, a= d=0\; mod \;{\mathfrak M}.$
In this case $\pi$ is of level less or equal to $r-1.$

\item ($C_1'$) $\begin{pmatrix}a&0\\0&d
\end{pmatrix}, a,d\in {\mathbb O} _{l'}, a= d\; mod \;{\mathfrak M}$ and $a\not=0 \;mod\; {\mathfrak M}.$
In this case $\pi=\tilde{\psi}_{\bar{a}}\otimes \pi'$ where $\pi'$ is a representation of $G^{(r)}$ of level less than $r-1$ and $\tilde{\psi}_{\bar{a}}$ is a primitive character of $G^{(r)}.$ In this case $\pi$ is of level  $r$ and of sublevel less or equal to $r-1.$

\item $(C_2$) $\begin{pmatrix}a&0\\0&d
\end{pmatrix}, a,d\in {\mathbb O} _{l'}, a\not= d\; mod \;{\mathfrak M}.$

\item ($C_3$) $\begin{pmatrix}0&1\\-\Delta&s
\end{pmatrix}, \Delta,s\in {\mathbb O} _{l'},  x^2-sx+\Delta$\; is irreducible  $mod \;{\mathfrak M}.$

\item $(C_4$) $\begin{pmatrix}a&1+b\\c&d
\end{pmatrix}, a, b, c, d \in{\mathbb  O}_{l'},  b, c, a-d\in {\mathfrak M}.$
Then $\pi=\tilde{\psi}_{\bar{a}}\otimes \pi'$ where $\pi'$ is a primitive representations which restriction $\pi'\vert_{K_l^{(r)}}=e\bigoplus_{\beta'\in \Omega'} \psi_{\beta'}$ where  
$\beta'$ is conjugated under $G^{(r)}$ to  $\beta_{C_4'(\Delta,s)}=\begin{pmatrix}0&1\\-\Delta&s
\end{pmatrix} \Delta,s\in {\mathbb O} _{l'}, \Delta,s\in {\mathfrak M}.$
\end{itemize}
In the case $C_2, C_3, C_4$ the representation is strongly primitive.

The orbit associated to $(C_1)$  and $(C_1')$ are not regular in the sense of Hill \cite{Hi1}. The orbits of type $(C_2), (C_3), (C_4)$ are regular. The orbit $(C_2)$ and $(C_4)$ are {\bf split} \cite{Hi1} whereas $(C_3)$ is {\bf cuspidal} in Hill's terminology \cite{Hi2}. In all the cases where $\beta$ is regular, we have $T(\psi_{\beta})=({\mathcal O}_{r}[\hat{\beta}])^\times K_{l'}$ where $\hat{\beta}$ is an element in $M_2({\mathcal O}_r)$ having $\beta$ as projection in $M_2({\mathcal O}_{l'}).$ 
\end{proposition}

In \cite{St}  a complete classification (valid even for $p=2$)  of irreducible representation of $GL(2,{\mathcal O}_r)$ is given using an inductive process. These representations fall in two classes: they can be strongly primitive of level $r$ or they are twisted by a character from an irreducible representation of level less or equal to $r-1$.
Therefore the knowledge of all irreducible strongly primitive representations of $GL(2,{\mathcal O}_k)$ for every $k\leq r$ gives, after twisting by characters, the complete list of irreducible representations of $GL(2,{\mathcal O}_r).$

We now proceed and study in detail the strongly primitive representations of $GL(2,{\mathcal O}_r)$.
They fall into classes according to the previous proposition.

\begin{defi}
The representations associated  to the orbits $C_2$ will be called {\bf  principal split} representations. 

The representations associated to the orbits $C_3$ will be called {\bf cuspidal } representations

The representations associated to the orbits $C_4$ will be called {\bf non-principal split }representations. 
\end{defi}

At this point we have to distinguish two cases:
\begin{itemize}
\item the simplest case  is when $r$ is even, i.e $l=l'=\frac{r}{2},$  in this case $K_{l'}$ is abelian.
\item the more complicated case is when   $r$ is odd,  i.e  $l'=l-1=\frac{r-1}{2},$ in this case $K_{l'}$ is not abelian.

\end{itemize}
{\bf In the rest of this work   we study the case where $r$ is even}, the case where $r$ is odd is studied in \cite{BaRo}.

We use the notations and results of of \cite{St}. Let $\beta$ be an element of $M_2({\mathcal O}_{l'})$ belonging to the orbits $C_2, C_3$ or $C_4$ and let $\hat{\beta}$ be any lift of $\beta$ in $M_2({\mathcal O}_r)$.
Let $\theta\in Hom({\mathcal  O}_r[\hat{\beta}]^\times, {\mathbb C}^\times),$ be a character such that  $\theta$ and $\psi_\beta$ coincide on 
$K_{l'}\cap {\mathcal  O}_r[\hat{\beta}]^\times$, one define $\theta\psi_\beta$ to be the one dimensional representation of $T(\psi_\beta)={\mathcal O}_{r}[\hat{\beta}]^\times K_{l'}$  by 
$(\theta\psi_\beta)(xy)=\theta(x)\psi_\beta(y), x\in {\mathcal  O}_r[\hat{\beta}]^\times, y\in K_{l'}.$ Then the representation $\pi(\theta,\beta)=Ind_{{\mathcal O}_{r}[\hat{\beta}]^\times K_{l'}}^{G^{(r)}}(\theta\psi_\beta)$ is an irreducible representation. It is shown that the set of representations $\pi(\theta, \beta)$ up to isomorphism depends only on the orbit of $\beta$ and is independent of the choice of lift $\hat{\beta}.$ Furthermore up to isomorphism this is the complete list of strongly primitive representations of level $r.$

We now proceed further and analyse in detail  these representations. We have tried to simplify as much the construction of these representations, this will be important for computing their characters. Note that there is a neat construction for the Principal split and Cuspidal representations, but the non principal split representations resist such a description..

\subsection{Principal Split representations.}

Let $a,d\in {\mathbb O}_{l},a\not=d\;\; \text{mod}\;{\mathfrak M},$ we define $\beta(C_2(a,d))=\begin{pmatrix}a&0\\0&d
\end{pmatrix}.$
Let $g=\begin{pmatrix}1+\varpi^l x&\varpi^l y\\ \varpi^l z&1+\varpi^l t
\end{pmatrix}\in K_l.$
We have $\psi_{\beta(C_2(a,d))}(g)=\psi(\varpi^l ax+\varpi^l t d).$
The number of characters of the form $\psi_{\beta(C_2(a,d))}$ is $\vert {\mathbb O}_l\vert \vert {\mathbb O}_l\setminus {\mathfrak M}\vert =q^l(q^l-q^{l-1})=(q-1)q^{r-1}.$
Because $\beta(C_2(a,d))$ and  $\beta(C_2(d,a))$ are the only elements in  the same orbit under the conjugation action,   the number of orbits of type $C_2$ is $\frac{1}{2}(q-1)q^{r-1}.$

We have $T(\psi_{\beta(C_2(a,d))})=S K_{l}$ with 
$S=\{\begin{pmatrix}s_1&0\\0& s_2\end{pmatrix}, s_1,s_2\in{\mathcal O}_r^{\times}\}=({\mathcal O}_r[\hat{\beta}])^{\times}$ where $\hat{\beta}$ is any lift of $\beta(C_2(a,d))$ in $M_2({\mathcal O}_r).$
Note that $T(\psi_{\beta(C_2(a,d))})=\{\begin{pmatrix}s_1&\varpi^l y\\ \varpi^l z& s_2\end{pmatrix}, s_1,s_2\in {\mathcal O}_r^{\times}, y,z\in {\mathcal O}_r\}=T(C_2)$ and is independent of $a, d.$

We have $\vert T(\psi_{\beta(C_2(a,d))})\vert=\vert{\mathcal O}_r^{\times}\vert ^2(q^l)^2=(q-1)^2q^{3r-2}.$

Because $K_{l}\cap S=\{\begin{pmatrix}1+\varpi^l x&0\\0& 1+\varpi^l y\end{pmatrix}, x,y \in {\mathcal O_r}\}$ we have $\vert K_{l}\cap S\vert=q^{l}q^{l}=q^r.$ 
The number of characters $\theta:S\rightarrow {\mathbb C}^\times$ which are equal to $\psi_{\beta(C_2(a,d))}$ on $K_{l}\cap S$ is given by $\frac{\vert S\vert}{\vert S\cap K_{l}\vert }=q^{r-2}(q-1)^2.$
The  irreducible principal split representation $\pi(\theta,\beta((C_2(a,d)))$ is   strongly primitive and of dimension 
$\frac{\vert G^{(r)}\vert }{\vert S K_{l}\vert}=(q+1) q^{r-1}.$

Moreover from the counting above, the  number of  inequivalent irreducible  principal split representations of $G^{(r)}$   is $\frac{1}{2}(q-1)^3q^{2r-3}$. 

We now give a precise description of these characters $\theta.$

A one dimensional representation  $\theta$ of $S$ is necessarily equal to $\theta_{\mu,\mu'}$ where  $\mu,\mu'$ are characters of ${\mathcal O}_r^\times$ and 
$\theta_{\mu,\mu'}(\begin{pmatrix}s_1&0\\0& s_2\end{pmatrix})=\mu(s_1)\mu'(s_2).$
The condition that $\theta$ and $\psi_{\beta(C_2(a,d))}$ are equal on $S\cap K_l$ is given by:
$$\mu(1+\varpi^l x)\mu'(1+\varpi^l y)=\psi(\varpi^l ax +\varpi^l d y), \forall x, y \in {\mathcal O_r},$$
which is equivalent to 
$$\mu(1+\varpi^l x)=\psi(\varpi^l ax), \forall x \in {\mathcal O_r}, \mu'(1+\varpi^l y)=\psi(\varpi^l d y), \forall y \in {\mathcal O_r}.$$
Note that we  can recover $a, d\in {\mathbb O_l}$ from the characters $\mu, \mu'$ using  the last condition.

In particular one has the important property that $\mu {\mu '}^{-1}(1+\varpi^l x)=\psi(\varpi^l x (a-d)), \forall x \in{\mathcal O_r}$. Asking that the  restriction of $\mu {\mu '}^{-1}$ to the multiplicative group 
$(1+\varpi^l {\mathcal O}_r)$ is of level $l$ is equivalent to the fact that $a\not=d\; mod\; {\mathfrak M}.$

We can therefore parametrize the set of principal split representations of $G^{(r)}$ by a pair of characters of ${\mathcal O}_r^\times.$ 
We say that  a couple $(\mu,\mu')$ of characters of ${\mathcal O}_r^\times$ is regular if and only if 
 $\mu {\mu '}^{-1}\vert_{1+\varpi^l {\mathcal O}_r}$ is primitive (i.e of level $l$).

 For such regular couple $(\mu,\mu'),$ we will denote $\Pi_{\mu,\mu'}$ the irreducible representation $\pi(\theta_{\mu,\mu'},\psi_{\beta(C_2(a,d))} )$ where
 $a \; (\text{resp.} \;d)$ are defined by $\mu \;(\text{resp.}\; \mu').$
The representation  $\Pi_{\mu,\mu'}$ and  $\Pi_{\mu',\mu}$ are isomorphic and up to equivalence depend only on the pair $\{\mu,\mu'\}$.

Remark: Note that the character $\theta_{\mu,\mu'}\psi_{\beta(C_2(a,d))}$ of $SK_{l}$, denoted $\mu\boxtimes\mu'$  has the following simple expression on $T(C_2)$
\begin{equation*}
(\mu\boxtimes\mu')(\begin{pmatrix}s_1&\varpi^l y\\ \varpi^l z& s_2\end{pmatrix})=\mu(s_1)\mu'(s_2), s_1,s_2\in {\mathcal O_r}^\times, x,y\in {\mathcal O_r},
\end{equation*}
and   $\Pi_{\mu,\mu'}=Ind_{T(C_2)}^{G^{(r)}}(\mu\boxtimes\mu').$ With this description we do not use $(a,d),$ which can be recovered from $(\mu,\mu').$

\subsection{Cuspidal representations.}

Let $F^{ur}$ be the maximal unramified extension of $F$, we have $Gal(F^{ur}/F)\simeq Gal(\bar{{\mathbb F}}_q/{\mathbb F}_q),$ let $\sigma$ be the element of $Gal(F^{ur}/F)$ corresponding to the Frobenius automorphism $Fr$ of $\bar{{\mathbb F}}_q.$ Let $E$ be the unique unramified extension of $F$ of degree $2$ i.e  $E=\{x\in F^{ur}, \sigma^2(x)=x\},$
we denote ${\mathcal O}^{E}$ the ring of integers  of $E, $ its maximal ideal is generated by $\varpi$ and its residual field is ${\mathbb F}_{q^2}.$

We denote ${\mathcal O}^{E}/\varpi^k  {\mathcal O}^{E}={\mathcal O}^{E}_k$ for $k$ positive integer.  We fix $r$ integer and denote ${\mathbb O}_k^E,$  for $0\leq k\leq r$ the image of compatible 
sections  of  ${\mathcal O}^{E}_k$. As usual we define the maps $Tr, N:E\rightarrow F$ by $Tr(x)=x+\sigma(x), N(x)=x\sigma(x).$

For  any $\tau \in {\mathbb O}^{E}_l$ such that $\tau-\sigma( \tau)\not=0$ mod  $\varpi$, we define 
$\beta(C_3(\tau))= \begin{pmatrix}
0&1\\-N(\tau)&Tr(\tau)
\end{pmatrix},$ which is a matrix of type $C_3.$
The reduction of $\beta(C_3(\tau))$ in $M_2(\Bbbk)$  is  $\overline {\beta(C_3(\tau))}=\begin{pmatrix}
0&1\\-\bar{\tau} Fr(\bar{\tau})&\bar{\tau}+Fr(\bar{\tau})
\end{pmatrix}$ with $Fr(\bar{\tau})\not=\bar{\tau}, $ i.e $\bar{\tau}\in {\mathbb F}_{q^2}\setminus {\mathbb F}_q.$

Let $g=\begin{pmatrix}1+\varpi^l x&\varpi^l y\\ \varpi^l z&1+\varpi^l t
\end{pmatrix}\in K_l,$
we have $\psi_{\beta(C_3(\tau))}(g)=\psi(\varpi^l z-\varpi^lN(\tau) y +\varpi^l Tr(\tau)t).$

The number of characters of the form $\psi_{\beta(C_3(\tau))}$  is $\frac{1}{2}(q-1)q^{r-1}.$
This is because the set $\{\tau \in {\mathbb O}^{E}_l ,\tau=\sigma( \tau) \;\text{mod } \;\varpi \}$ is of cardinal $qq^{2(l-1)}$ and  
$\psi_{\beta(C_3(\tau))}=\psi_{\beta(C_3(\sigma(\tau)))}.$ Therefore the number of orbits of type $C_3$ is $\frac{1}{2}(q-1)q^{r-1}.$

Let $\hat{\tau}$ a representative of $\tau$ in  ${\mathcal O}^{E}_r,$
 we define  $\hat{\beta}$ a lift of $\beta(C_3(\tau))$, $\hat{\beta} =\begin{pmatrix}
0&1\\-\hat{\tau} \sigma(\hat{\tau})&\hat{\tau}+\sigma(\hat{\tau})
\end{pmatrix}.$

$T(\psi_{\beta(C_3(\tau))})=S K_l$ and we have 
$S={\mathcal O}_r[\hat{\beta}]^{\times }=\{\begin{pmatrix} a&b\\-b\hat{\tau} \sigma(\hat{\tau})&a+b(\hat{\tau}+\sigma(\hat{\tau}))
 \end{pmatrix}, a,b\in {\mathcal O}_r, (a+b\hat{\tau}) (a+b\sigma(\hat{\tau}))\in {\mathcal O}_r^{\times}
\}$.
As a result we obtain $$T=\{\begin{pmatrix} x&y\\-y \hat{\tau} \sigma(\hat{\tau})+\varpi^l z&x+y(\hat{\tau}+\sigma(\hat{\tau}))+\varpi^l t
 \end{pmatrix}, x,y,z,t\in {\mathcal O}_r \}.$$

Let $a, b\in {\mathcal O}_r, $ because $\bar{\tau}\in {\mathbb F}_{q^2}\setminus {\mathbb F}_q$, we have the equivalence: $(a+b\hat{\tau}) (a+b\sigma(\hat{\tau}))=0\;mod\; \varpi$ if and only if $a$ and $b$ are equal to $0$ mod $\varpi$.

Therefore 
$S=\{\begin{pmatrix} a&b\\-b\hat{\tau} \sigma(\hat{\tau})&a+b(\hat{\tau}+\sigma(\hat{\tau}))
 \end{pmatrix}, a,b\in {\mathcal O}_r \setminus \varpi {\mathcal O}_r \}$ implying $\vert S\vert  =(q^2-1)q^{2(r-1)}.$

Moreover $S\cap K_l=\{\begin{pmatrix}
a&b\\-b\hat{\tau} \sigma(\hat{\tau})&a+b(\hat{\tau}+\sigma(\hat{\tau}))
\end{pmatrix}, a,b\in {\mathcal O}_r\; a=1\; mod\; \varpi^l, \;b=0\; mod\; \varpi^l
\}$

We therefore have $\vert S\cap K_l \vert=(q^{r-l})^2=q^r, $ from which it follows that 
$\vert T(\psi_{\beta(C_3(\tau))})\vert=(q^2-1)q^{3r-2}.$

The number of characters $\theta:S\rightarrow {\mathbb C}^\times,$ which extend $\psi_{\beta(C_3(\tau))}$ on $S\cap K_l,$ is given by $\frac{\vert S\vert}{\vert S\cap K_{l}\vert }=(q^2-1)q^{r-2}.$

From the counting argument above, the number of inequivalent cuspidal  representations of $G^{(r)}$ is $\frac{1}{2}(q-1)(q^2-1)q^{2r-3}$. 
These representations  are all strongly primitive and of dimension 
$\frac{\vert G_r\vert }{\vert S K_{l}\vert}=(q-1) q^{r-1}.$

We now give a precise description of these characters $\theta.$

Let $\nu,\nu'$ characters $({\mathcal O}_r^E)^{\times}\rightarrow {\mathbb C}^\times,$ we define  a character
 $\theta_{\nu,\nu'}:S\rightarrow {\mathbb C}^\times$ by 
\begin{equation*}
\theta_{\nu,\nu'}(\begin{pmatrix} a&b\\-b\hat{\tau} \sigma(\hat{\tau})&a+b(\hat{\tau}+\sigma(\hat{\tau}))
 \end{pmatrix})=\nu(a+b \hat{\tau})\nu'(a+b \sigma(\hat{\tau})).
\end{equation*}
Each character of $S$ is of this type.

The condition that $\theta$ and $\psi_{\beta(C_3(\tau))}$ are equal on $S\cap K_l$ is given by:
\begin{equation}
\label{conditionC3even}\nu(1+\varpi^l x+\varpi^l y\tau)\nu'(1+\varpi^l x+\varpi^l y\sigma(\tau))=
\psi(\varpi^l {\tau }(x+y{\tau})+\varpi^l \sigma ({\tau}) (x+y\sigma({\tau}))), \forall x,y\in {\mathcal O}_r.
\end{equation}

Let ${\mathcal O}_r^{E(\pm)}= \{z\in {\mathcal O}_r^{E}, \sigma(z)=\pm z\},$
 $(1+\varpi^l{\mathcal O}_r^{E(\pm)},\times)$ are  subgroups of $({\mathcal O}_{r}^E)^{\times}.$ 
The  condition (\ref{conditionC3even}) implies that 
\begin{equation*}\nu(1+\varpi^l z )\nu'(1+\epsilon\varpi^l z)=
\psi(\varpi^l z({\tau }+\epsilon \sigma({\tau })), \forall z\in {\mathcal O}_r^{E(\pm)},\forall \epsilon\in \{+,-\}.
\end{equation*}
In particular we obtain that $\nu{\nu'}^{-1}(1+\varpi^l z)=\psi(\varpi^l(\tau-\sigma(\tau))z)$ for $z\in {\mathcal O}_{r}^{E(-)},$ i.e $\nu{\nu'}^{-1}$ is a representation of $(1+\varpi^l{\mathcal O}_r^{E(-)}) $ of level $l.$ 

We say that a couple $(\nu,\nu')$ of characters of $({\mathcal O}_r^E)^{\times}$ is regular if and only if $\nu{\nu'}^{-1}$ is a representation of $ (1+\varpi^l{\mathcal O}_r^{E(-)})$ of level $l.$ 

In the case where $p\not=2,$  we can recover $\tau\in {\mathbb O}^{E}_l$ from the knowledge of $\nu,\nu'.$

Indeed if $(\nu,\nu')$ is a regular couple of  characters of $({\mathcal O}_r^E)^{\times}, $ let $\epsilon\in \{+,-\}$ then there exists a unique $\tau^{(\epsilon)}\in {\mathbb O}_l^E$ such that 
$\nu(1+\varpi^l z )\nu'(1+\epsilon\varpi^l z)=\psi(\varpi^l \tau^{(\epsilon)}),  \forall z \in {\mathcal O}_r^{E(\epsilon)} .$  The regularity condition 
   implies that $\tau^{(-)}\notin {\mathfrak M}.$
It is an easy exercise to show that if we define $\tau=\frac{1}{2}(\tau^{(+)}+\tau^{(-)})$ then the condition (\ref{conditionC3even}) holds.

Given $(\nu,\nu')$  a regular couple of  characters of $({\mathcal O}_r^E)^{\times}$,  we will denote  $\nu\boxtimes \nu'= 
 \theta_{\nu,\nu'} \psi_{\beta(C_3(\tau))} ,$  where $\tau$ is defined from $(\nu,\nu').$ 

 We denote 
${\mathcal C}_{\nu,\nu'}$ the representation 
$ \pi(\theta_{\nu,\nu'}, \psi_{\beta(C_3(\tau))})=Ind_{T(C_3(\tau))}^{G^{(r)}}(\nu\boxtimes \nu').$   Note that ${\mathcal C}_{\nu,\nu'}$ is isomorphic to  
${\mathcal C}_{\nu'\circ\sigma,\nu\circ\sigma}.$

Remark: At this point, we want to make a connection with the work of \cite{LCW}. They are working with the case where $p\not=2$ and ${\mathcal O}={\mathbb Z}_p$. They have chosen a different   representative of $C_3$ and a different choice of $\beta$ associated to $C_3.$
For $F$ any local field with $p\not=2$, denote ${\mathcal O}$ the ring of integer of $F,$
if  $\tilde{\rho}\in {\mathbb O}_l$ and $\tilde{\epsilon}\in{\mathbb O}_l$ such that $\tilde{\epsilon}$ is not a square in ${\mathcal O}_r$ and is invertible in ${\mathcal O_r}$ we define 
 $\tilde{\beta}=\begin{pmatrix}\tilde{\rho}&\tilde{\epsilon}\\1&\tilde{\rho}\end{pmatrix}.$
(The elements $\tilde{\rho},\tilde{\epsilon}$ are the element $\alpha,\epsilon$  (page 1297) of \cite{LCW} in the case where ${\mathcal O}={\mathbb Z}_p).$
We have $\tilde{\beta}=P\psi_{\beta(C_3(\tau))}P^{-1}$ with $P=\begin{pmatrix}-\tilde{\rho}&1\\1&0\end{pmatrix},$ with $\tilde{\rho}=\frac{1}{2}(\tau+\sigma(\tau)), \tilde{\epsilon}=(\frac{1}{2}(\tau-\sigma(\tau)))^2.$
Note  they have used the same $\epsilon$ to parametrize the conjugacy classes of type $C_3$ as well as   the representions of cuspidal type.
We will prefer to proceed as follows: once forall, we have fixed  an invertible element $\epsilon$ in ${\mathcal O}_r$ such that $\epsilon$ is not a square. This $\epsilon$ is used for the proposition (\ref{theoclassificationconjugyclassthree}).  We will denote  $\Phi$ a square root in ${\mathcal O}_r^{E}$ of $\epsilon.$ 
The cuspidal representations are labelled by a regular couple $(\nu,\nu')$. This couple   defines $\tau\in {\mathbb O}_l^E$  from which we define  
$\tilde{\epsilon}=(\frac{1}{2}(\tau-\sigma(\tau)))^2$: $\epsilon$ is fixed whereas $\tilde{\epsilon}$ depends on the choice of the regular couple $(\nu,\nu').$
$\epsilon\tilde{\epsilon}$ is a square in ${\mathcal O}_r$  and we have $\epsilon\tilde{\epsilon}=u^2$ with $u=\frac{1}{2}\Phi(\tau-\sigma(\tau))\in {\mathcal O}_r^\times.$

\subsection{Non-Principal Split representations.}

Let $\Delta, s\in  {\mathbb O}_l\cap {\mathfrak M},$ we define $\beta_{C_4'(\Delta,s)}=\begin{pmatrix}0&1&\\-\Delta&s
\end{pmatrix}.$

Let $g=\begin{pmatrix}1+\varpi^l x&\varpi^l y\\ \varpi^l z&1+\varpi^l t
\end{pmatrix}\in K_l,$
we have $\psi_{\beta(C_4'(\Delta,s))}(g)=\psi(\varpi^l z-\Delta \varpi^l y +s\varpi^l t).$

The number of characters of the form $\psi_{\beta(C_4'(\Delta,s))}$  is $(q^{(r-l-1)})^2=q^{r-2}.$

Let $\hat{\Delta},\hat{s}$ lifts of $\Delta, s$ in ${\mathcal O}_r$ and define 
$\hat{\beta}(\hat{\Delta},\hat{s})=
\begin{pmatrix}0&1\\-\hat{\Delta}&\hat{s}
\end{pmatrix}$ 
a lift of $\beta_{C_4'(\Delta,s)}.$ 

We have $T(\psi_{\beta(C_4'(\Delta,s))})=SK_{l}$ with

 $S(\hat{\Delta},\hat{s})=S=({\mathcal O}_r[\hat{\beta}])^\times=
\{\begin{pmatrix}a&b\\
-\hat{\Delta} b& a+\hat{s} b\end{pmatrix}, a, b\in {\mathcal O}_r, a^2+\hat{s}ab+\hat{\Delta} b^2\in {\mathcal O}_r^{\times} \}$=
$\{\begin{pmatrix}a&b\\
-\hat{\Delta} b& a+\hat{s} b\end{pmatrix},  a\in {\mathcal O}_r^{\times}, b\in {\mathcal O}_r\}$

We have $\vert ({\mathcal O}_r[\hat{\beta}])^{\times}  \vert =
\vert {\mathcal O}_r\vert \vert {\mathcal O}_r\setminus \varpi {\mathcal O}_r\vert =(q-1)q^{2r-1}.$

Because $S\cap K_{l}=\{ \begin{pmatrix}a&b\\
-\hat{\Delta} b& a+\hat{s} b\end{pmatrix}, a, b\in \varpi^l {\mathcal O}_r\},$ $\vert S\cap K_l\vert= q^r$  therefore 
$\vert T(\psi_{\beta(C_4'(\Delta,s))}) \vert=\frac{\vert S\vert\vert K_l\vert}{\vert S\cap K_l \vert}=(q-1)q^{3r-1}.$

The number of characters $\theta:S\rightarrow {\mathbb C}^\times$ which extend $\psi_{\hat{\beta}}$ on $S\cap K_l$ is given by $\frac{\vert S\vert}{\vert S\cap K_{l}\vert }=(q-1)q^{r-1}.$

Finally the number of inequivalent irreducible representations of type $C_4'$ is $(q-1)q^{r-1}\times q^{r-2}=(q-1)q^{2r-3}$. These representations are all strongly primitive and of dimension 
$\frac{\vert G_r\vert }{\vert S K_{l}\vert}=(q^2-1) q^{r-2}.$

Finally these representations can be tensored with one dimensional characters of $G^{(r)}$ of the type ${\tilde\psi}_{a}$, $a\in {\Bbbk}$,  which give all the inequivalent representations of type $C_4.$ There is therefore $(q-1)q^{2r-2}$ inequivalent representations of this type, all of them being strongly primitive and of dimension $(q^2-1)q^{r-2}.$

We now give a precise description of the characters $\theta$ extending  $\psi_{\hat{\beta}}$  on  $S\cap K_l$ .
Because $\hat{\beta}^2=-\hat{\Delta} I+\hat{s}\hat{\beta}$ the group law on $S$ is given by 
$$(aI+b\hat{\beta})(a'I+b'\hat{\beta})=(aa'-\hat{\Delta} b b')I+(ab'+a' b+\hat{s}b b')\hat{\beta}$$ with $a,a'\in {\mathcal O}_r^{\times}, b,b'\in {\mathcal O}_r.$
Let $\theta:S\rightarrow {\mathbb C}^\times$ be a character, the center of $S$ being $Z(S)=\{aI,a\in{\mathcal O}_r^\times \},$ the restriction $\theta\vert_{Z(S)}$ defines a multiplicative character, denoted $\sigma$ of ${\mathcal O}_r^\times.$
Therefore we have $\theta(aI+b\hat{\beta})=\sigma(a)\eta(b/a)$ where $\eta:{\mathcal O}_r\rightarrow  {\mathbb C}^\times$ defined by $\eta(x)=\theta(1+x\hat{\beta}),x\in  {\mathcal O}_r$ satisfies:
\begin{equation*}
\eta(x)\eta(y)=\sigma(1-\hat{\Delta}xy)\eta(x\star y)
\end{equation*}
with 
\begin{equation}\label{staraddition}
x\star y=\frac{x+y+\hat{s}xy}{1-\hat{\Delta} xy}.\end{equation}

$({\mathcal O}_r, \star)$ is a commutative group.

Let $c:{\mathcal O}_r\times {\mathcal O}_r\rightarrow {\mathbb C}^\times$, be the map defined by $c(x,y)=\sigma(1-\hat{\Delta} x y),$
$c$ is a two cocycle in the sense that $c(x\star y,z)c(x,y)=c(x, y\star z)c(y,z), \forall x, y, z\in {\mathcal O}_r,$ and $\eta$ is therefore a projective representation of the additive group $({\mathcal O}_r,\star)$ associated to the $2$-cocycle $c^{-1}.$
The condition that $\theta$ extends $\psi_{\hat{\beta}}$ reads:

\begin{eqnarray*}\label{sigma}
&&\theta((1+\varpi^l a)I+\varpi^l b \hat{\beta})=\psi (-2\Delta \varpi^l b +s\varpi^l a+s^2\varpi^l b)\;\;\forall a,b\in {\mathcal O}_r\\
\label{eta}&&=\sigma(1+\varpi^l a )\eta(\varpi^l b/(1+\varpi^l a))=\sigma(1+\varpi^l a )\eta(\varpi^l b)
\end{eqnarray*}
which is equivalent to 
\begin{eqnarray*}
&\sigma(1+\varpi^l a)=\psi(s\varpi^l a),\forall  a\in {\mathcal O}_r\\
&\eta(\varpi^l b)=\psi((s^2-2\Delta)\varpi^l b),\forall  b\in {\mathcal O}_r.
\end{eqnarray*}
Therefore the restriction of $\sigma$ to the  group $1+\varpi^l {\mathcal O}_r$ is not primitive and it defines uniquely $s\in  {\mathbb O}_l\cap {\mathfrak M}.$
Note that the restriction of $\star$ and $+$ to $\varpi^l {\mathcal O}_r$ coincide  and the restriction of the cocycle $c$ to $\varpi^l {\mathcal O}_r$ is trivial. Therefore the restriction of $\eta$ to $(\varpi^l {\mathcal O}_r, +)$ is a one dimensional representation. It  is not primitive, and once $s$ is known through the knowledge of $\sigma,$ it defines uniquely $\Delta\in {\mathbb O}_l\cap {\mathfrak M}$ when $p\not=2.$

Remark 1.  One can endow ${\mathcal O}$ with an internal law $\star$ defined by the same formula as (\ref{staraddition}), but with $\hat{s},\hat{\Delta}\in \varpi{\mathcal O}$. By the theory of formal group it is shown that $({\mathcal O},\star)$ is isomorphic to  $({\mathcal O},+)$ when $F$ is of $0$ characteristic (Exercice 2 of \cite{Ne} page 345). This is  however not the case when we consider ${\mathcal O}_r$ and this prevent us to construct all the characters of $S$ using only projective characters of $({\mathcal O}_r, +).$

Remark 2. In \cite{LCW} it is said in the introduction that their methods could be applied to find the character value of $GL(2,{\mathcal O})$ which is the content of  our work. They say that it is easier for them to count the number of irreducible representations in the case where ${\mathcal O}={\mathbb Z}_p$ but from the work of \cite{St} there is no counting argument involved because the list of irreducible representations is known to be complete by Clifford theory.
Nevertheless we can check by a counting argument that the list of representations is indeed complete as follows.
The number of irreducible representations that we have constructed which fall into the classes of principal or non principal split and cuspidal  representations is 
$b''_r=\frac{1}{2}(q-1)^3q^{2r-3}+\frac{1}{2}(q-1)(q^2-1)q^{2r-3}+(q-1)q^{2r-2}=(q-1)q^{2r-1}$ which is equal to $n_r-qn_{r-1}=b'_r.$ Therefore the list of strongly primitive representations is complete.

\medskip
\medskip

 In the case where $p\not=2$ we can give a somewhat simpler descriptions of nonprincipal split representations  which is closer to the classification given in \cite{LCW}.
For $\tilde{\Delta},\tilde{s}\in {\mathbb O}_l\cap {\mathfrak M},$  denote $\beta_{C''_4(\tilde{\Delta},\tilde{s})}=\begin{pmatrix}\tilde{s}/2&1\\-\tilde{\Delta}&\tilde{s}/2\end{pmatrix}.$
The matrix $\beta_{C'_4(\Delta,s)}$ is conjugated under $G^{(r)}$ to the matrix $\beta_{C''_4(\Delta-\frac{s^2}{4},s)}.$ 
\begin{proposition}
We use the same notation as in proposition \ref{classificationorbitbeta}.
If $\pi$ is an irreducible representation of $G^{(r)}$  which orbit has a representative $\beta_{C'_4(\Delta,s)},$  with $\Delta,s\in {\mathbb O}_l\cap{\mathfrak M}$, we have   $\pi=\tilde{\psi}_{s/2}\otimes \pi''$ where $\pi''$ is a strongly primitive representation which restriction $\pi''\vert_{K_l^{(r)}}=e\bigoplus_{\beta''\in \Omega''} \psi_{\beta''}$ where  
$\beta''$ is conjugated under $G^{(r)}$ to  $\beta_{C_4''(\Delta-\frac{s^2}{4},0)}$ and $\tilde{\psi}_{\frac{s}{2}}$ is a  one dimensional representation of $G^{(r)}$ of the form $ \tilde{\psi}_{\frac{s}{2}}=\tilde{\chi}\circ det$  with $\tilde{\chi}$ one dimensional representation of ${\mathcal O}_r^{\times}$ and the restriction of $\tilde{\psi}_{\frac{s}{2}}$ to $K_l$ is given by $\psi_{\frac{s}{2}I}.$
\end{proposition}
\proof  Same proof as in \cite{St}.
\endproof

Therefore we have "absorbed"  $s$ by tensoring with a one dimensional representation, and we are left with the analysis of the representations associated to  $\beta_{C_4''(\Delta,0)}.$
Let $\hat{\Delta}$ be a lift in $\varpi{\mathcal O_r}$ of $\Delta.$
Let $\theta$ be a character of $S(\hat{\Delta},0)\rightarrow {\mathbb C}^\times.$ One can associate to it $\eta,\sigma$ satisfying the same relations as (\ref{sigma}, \ref{eta}) but with $\hat{s}=0$. In particular $\sigma$ is trivial on $1+\varpi^l {\mathcal O}_r.$
The representation $\pi(\theta, \beta_{C_4''(\Delta,0)})$ will be denoted ${\Xi}_{\Delta,\theta}.$  Up to isomorphism it does not depend on the choice of the lift of ${\Delta}.$
In order to obtain the complete set of representations which are non principal split, we have to tensor them with $\tilde{\psi}_{\frac{s}{2}}$ $s\in {\mathbb O}_l\cap {\mathfrak M}$ and with $\tilde{\psi}_{\bar{a}_0}, a_0\in {\mathbb O}_1.$

We shall also denote ${\Xi}_{a,\Delta,\theta}=\tilde {\psi}_{a}\otimes {\Xi}_{\Delta,\theta},$   with $a\in {\mathbb O}_l$ and $\tilde{\psi}_{a}$ is a one dimensional representation of $G^{(r)}$ of the type $\tilde{\psi}_{a}=\tilde{\chi}_{a}\circ det$ where $\tilde{\chi}_{a}$ a one dimensional representation of ${\mathcal O}_r^\times$ and the restriction 
of $\tilde{\psi}_{a}$ to $K_l$ is given by $\psi_{aI}$, we have $a=a_0+s/2$ with $a_0\in {\mathbb O}_1$ and $s\in {\mathbb O_l}\cap {\mathfrak M}.$
\bigskip

The following table summarizes the essential informations on strongly primitive representations:

\medskip

\begin{tabular}{|c|c|c|}\hline
Strongly primitive irrep of odd level $r$  & Dimension  &  Number of inequivalent irrep \\
\hline\hline
Principal Split representations  $\Pi_{\mu,\mu'}$ & $(q+1)q^{r-1}$& $\frac{1}{2}(q-1)^3q^{2r-3}$ \\
\hline Cuspidal representations ${\mathcal C}_{\nu,\nu'}$& $(q-1)q^{r-1}$ & $\frac{1}{2}(q-1)(q^2-1)q^{2r-3}$ \\
\hline
Non Principal Split representations  ${\Xi}_{a,\Delta,\theta}$ & $(q^2-1)q^{r-2}$& $(q-1)q^{2r-2}$\\
\hline
\end{tabular}
\bigskip

\section{Characters}
We will use the formula  of Frobenius giving the character of an induced representation. Let  $G$  a finite group, $H$ a subgroup of $G$ and $\pi$ a finite dimensional complex representation of $H$ having character $\chi_{\pi}$, then the character of $Ind_{H}^G(\pi)$ is given by:
\begin{eqnarray}
tr(g\vert_{Ind_{H}^G(\pi)})&=&\frac{1}{\vert H \vert }\sum_{t\in G}\chi_{\pi}^{0}(tgt^{-1}), \;\forall g\in G\\
&=& \sum_{t\in X}\chi_{\pi}^{0}(tgt^{-1}), \end{eqnarray}
where for any function $\phi$ on $H$, $\phi^{0}$ denotes the extension of $\phi$ on $G$ by $\phi^{0}(g)=\phi (g)$ if $g\in H$ and zero otherwise, and $X$ denote any section of the right cosets of $H$ in $G$.
We will apply this formula to the case where $G=G^{(r)}$, $r$ even  and $H=T(\psi_{\beta})=SK_{l}.$ 

We will assume that $p\not=2$ in the rest of this work 
 and we will use the proposition (\ref{theoclassificationconjugyclassthree}) to obtain a representative set of conjugacy classes.

Remark: Although representatives of conjugacy classes are known in the case  $p=2$ as well as an exhaustive list of irreducible representations, computing characters using   Frobenius formula appears to be very complicated.

\subsection{ Principal split representations}
Having chosen $\beta$ ot type $C_2$, we obtain:
\begin{proposition}
A section of the right cosets of $SK_{l}$  is given by the following set of matrices
 $X=\{ e_{x,y}, f_{x,z}; x,y\in {\mathbb O}_{l},  z\in{\mathbb  O}_{l-1}\}$ where 
$ e_{x,y}=\begin{pmatrix}1&x\\0&1\end{pmatrix} \begin{pmatrix}1&0\\y&1\end{pmatrix} ,
f_{x,z}=\begin{pmatrix}1&x\\0&1\end{pmatrix} \begin{pmatrix}\varpi z&1\\1&\varpi z\end{pmatrix} .$
\end{proposition}
\proof 
We have $SK_{l}=\{\begin{pmatrix} s_1&\varpi^l y\\\varpi^l z& s_2\end{pmatrix}, s_1,s_2\in {\mathcal O}_r^{\times}, y,z\in { \mathcal O}_r
\}.$ It is an easy exercise to show that the right coset  associated to the elements of $X$ are disjoints. Moreover $\vert X\vert=(q^{l})^2+q^{l}q^{l-1}=\vert G/SK_l\vert$ therefore $X$ is a section of the right cosets of  $SK_{l}$.
\endproof
We will denote $\xi:{\mathcal O}_r\rightarrow {\mathcal O}_r^\times$ the function defined by $\xi(z)=1-\varpi^2 z^2=det(f_{x,z}).$
Frobenius formula therefore gives:
\begin{equation}
tr(g\vert_{\Pi_{\mu,\mu'}})=\sum_{t\in X}(\mu\boxtimes\mu')^0(tgt^{-1})
=S_{e}(g)+S_f(g),
\end{equation}
where $
S_{e}(g)=\sum_{x,y\in {\mathbb O_{l}}}(\mu\boxtimes\mu')^0(e_{x,y}ge_{x,y}^{-1})$ and $S_f(g)=\sum_{x\in {\mathbb O_{l}},z\in {{\mathbb  O}_{l-1}}}(\mu\boxtimes\mu')^0(f_{x,z}gf_{x,z}^{-1}).$
\medskip

$\bullet$ Conjugacy class of type $I.$

$I_{\alpha}$ being central we have $tr(I_\alpha\vert_{\Pi_{\mu,\mu'}})=\frac{\vert G\vert}{\vert SK_{l}\vert}\mu(\alpha)\mu'(\alpha)=q^{r-1}(q+1)\mu(\alpha)\mu'(\alpha).$

\medskip

$\bullet$ Conjugacy class of type ${ C}.$
\begin{proposition}
$tr(C_{i,\alpha,\beta}\vert_{\Pi_{\mu,\mu'}})=0.$
\end{proposition}
\proof
$e_{x,y}C_{i,\alpha,\beta}e_{x,y}^{-1}=\begin{pmatrix}
\alpha-\varpi^i \epsilon\beta y +\varpi^i\beta x (1-\epsilon y^2)& (1+xy)^2\varpi^i\epsilon \beta-x^2\varpi^i\beta\\
\varpi^i\beta (1-\epsilon y^2)& \alpha+\varpi^i\epsilon\beta y-\varpi^i\beta x (1-\epsilon y^2)\end{pmatrix}.$
Therefore this matrix does not belong to $SK_l$ if $i<l$ (because $\beta$ is invertible) and  its value on $(\mu\boxtimes\mu')^0$ is $0.$ 
When $i\geq l$, we necessarily  have $\alpha\in {\mathcal O}_r^{\times},$
\begin{eqnarray*}
S_e(C_{i,\alpha,\beta})&=&\sum_{x,y\in {\mathbb O_{l}}}
\mu(\alpha-\varpi^i \epsilon\beta y+\varpi^i\beta x(1-\epsilon y^2))
\mu'(\alpha+\varpi^i \epsilon\beta y-\varpi^i\beta x(1-\epsilon y^2))\\
&=&\mu(\alpha)\mu'(\alpha)\sum_{x,y\in {\mathbb O_{l}}}
(\mu\mu'{}^{-1})(1-\frac{\varpi^i}{\alpha} (\epsilon\beta y+\beta x(1-\epsilon y^2))).
\end{eqnarray*}
Because the  restriction of $\mu {\mu '}^{-1}$ to the multiplicative group 
$(1+\varpi^l {\mathcal O}_r)$ is primitive and $i<r$,  we have $\sum_{x\in {\mathbb O_{l}}}
(\mu\mu'{}^{-1})(1-\frac{\varpi^i}{\alpha} (\epsilon\beta y+\beta x(1-\epsilon y^2)))=0,
$
therefore $S_e(C_{i,\alpha,\beta})=0.$
We have 
\begin{eqnarray*}&&
f_{x,z}C_{i,\alpha,\beta}f_{x,z}^{-1}=\\
&&
\begin{pmatrix}\alpha+\frac{\varpi^i}{\xi(z)}(-\beta(1+\varpi zx)\varpi z +\epsilon\beta (\varpi z +x))& \star\\
\frac{\varpi^i}{\xi(z)}\beta (\epsilon-\varpi^2 z^2)& \alpha +\frac{\varpi^i}{\xi(z)}(\varpi z\beta (1+x\varpi z)-\epsilon\beta(x+\varpi z)),
\end{pmatrix}
\end{eqnarray*}
in order to have a non zero value by  $(\mu\boxtimes\mu')^0$ it is necessary that $i\geq l.$ In that case $\alpha$ is invertible and 

\begin{eqnarray*}
&&S_f(C_{i,\alpha,\beta})=\\
&&\sum_{x\in {\mathbb O_{l}}, z\in {\mathbb O_{l-1}}}
\mu(\alpha-\varpi^i \frac{(\beta(1+\varpi z x)\varpi z-\epsilon \beta (\varpi z+x))}{\xi(z)} )
\mu'(\alpha+\varpi^i \frac{(\beta(1+\varpi z x)\varpi z-\epsilon \beta (\varpi z+x))}{\xi(z)}) \\
&=&\mu(\alpha)\mu'(\alpha)\sum_{x\in {\mathbb O_{l}}, z\in {\mathbb O_{l-1}}}
(\mu\mu'{}^{-1})(1-\frac{\varpi^i}{\beta\alpha} (\frac{(\varpi z-\epsilon \beta \varpi z)+x(\varpi^2z^2-\epsilon\beta)}{\xi(z)})).
\end{eqnarray*}
The sum over $x$ gives again $0$ for the same reason as the one used for proving that   $S_e(C_{i,\alpha,\beta})=0.$
Therefore $S_f(C_{i,\alpha,\beta})=0.$ The evaluation of the character on elliptic elements is $0.$
\endproof

$\bullet $
Conjugacy class of type $D.$

\begin{lemma}\label{Identity1}
Let $\lambda: 1+\varpi^l {\mathcal O}_r\rightarrow {\mathbb C}^\times$ be a primitive  character with $r=2l.$ Let $u\in {\mathcal O}_r^{\times}$ and 
$i\in \llbracket 0,r-1\rrbracket, $ the following identity holds:
\begin{equation*}
\sum_{x,y\in{\mathbb O}_l \atop v(y)\geq l-i, v(x)\geq l-i}\lambda (1+u\varpi^i xy)=q^i.
\end{equation*}

\end{lemma}
\proof 
Up to changing the primitive character, we can always assume $u=1.$

If $i\geq l$ then there is no constraint on $x,y$. We can write $i=l+i'$ and we have to evaluate 
$A=\sum_{x,y\in{\mathbb O}_l }\lambda(1+\varpi^l x\varpi^{i'}y).$
$y$ being fixed, the sum over $x$ is nul unless $v(\varpi^{i'}y)\geq l$ and in this case it gives $\vert {\mathbb O}_l\vert.$
The summation over $x,y$ gives therefore
 $A=\vert {\mathbb O}_l\vert \vert \{y\in {\mathbb O}_l , v(\varpi^{i'}y)\geq l\}\vert=
\vert {\mathbb O}_l\vert \vert \varpi^{l-i'}{\mathbb O}_{i'}\vert=q^lq^{i'}=q^i.$

If $i<l$, we set $x=\varpi^{l-i}x', y=\varpi^{l-i}y', x',y'\in {\mathbb O}_i$, therefore 
$A=\sum_{x',y'\in {\mathbb O}_i}\lambda(1+ \varpi^{l}\varpi^{l-i}x' y').$
Noting that $z\mapsto \lambda(1+ \varpi^{l}\varpi^{l-i}z)$ is a primitive character of ${\mathcal O}_i$, we obtain that for $y'$ fixed the sum over $x'$ gives $0$ unless $y'=0.$  The summation over $x',y'$ gives therefore $A =\vert {\mathbb O}_i\vert=q^i.$

\endproof

\begin{proposition}
$tr(D_{i,\alpha,\delta}\vert_{\Pi_{\mu,\mu'}})=q^{i}(\mu(\alpha)\mu'(\delta)+ \mu'(\alpha) \mu(\delta) ).$
\end{proposition}

\proof
From 
\begin{equation*}
e_{x,y}D_{i,\alpha,\delta}e_{x,y}^{-1}=
\begin{pmatrix}
\alpha+(\alpha-\delta)xy&(\delta-\alpha)x(1+xy)\\
(\alpha-\delta)y&\delta-(\alpha-\delta)xy
\end{pmatrix},
\end{equation*}
we see that only  the matrix such that $v((\alpha-\delta)y)\geq l$  and $v((\delta-\alpha)x(1+xy))\geq l$   can  contribute to $S_e(D_{i,\alpha,\delta}).$
Because $v(\alpha-\delta)=i$, we necessarily have $v(y)\geq l-i.$
As a result we have to distinguish two cases: $i \geq l$ or $i< l.$

In the first case $i\geq l$, there is no condition on $x,y$.
In the second case we necessarily have $v(y)\geq l-i>0$, therefore $1+xy$ is invertible and we necessarily have $v(x)\geq l-i.$
Therefore:
\begin{eqnarray*}
S_e(D_{i,\alpha,\delta})&=&\sum_{x,y\in{\mathbb O}_l\atop v(y)\geq l-i, v(x)\geq l-i}\mu(\alpha+(\alpha-\delta)xy)\mu'(\delta-(\alpha-\delta)xy)\\
&=&  \mu(\alpha)\mu'(\delta)
\sum_{x,y\in{\mathbb O}_l\atop  v(y)\geq l-i, v(x)\geq l-i}\mu(1+\frac{(\alpha-\delta)}{\alpha}xy)\mu'(1-\frac{(\alpha-\delta)}{\delta}xy)\\
&&= \mu(\alpha)\mu'(\delta)\sum_{x,y\in{\mathbb O}_l \atop v(y)\geq l-i, v(x)\geq l-i}\mu(1+\frac{(\alpha-\delta)}{\alpha}xy)\mu'(1-\frac{(\alpha-\delta)}{\alpha}xy)
\end{eqnarray*}
the last equality holds because the equality  $\frac{(\alpha-\delta)}{\delta}xy=\frac{(\alpha-\delta)}{\alpha}xy$ holds from the valuation condition on $x,y.$
As a result:
\begin{eqnarray*}
S_e(D_{i,\alpha,\delta})&=&\mu(\alpha)\mu'(\delta)\sum_{x,y\in{\mathbb O}_l\atop  v(y)\geq l-i, v(x)\geq l-i}(\mu\mu'^{-1})(1+\frac{(\alpha-\delta)}{\alpha}xy)\\
&&=q^i\mu(\alpha)\mu'(\delta),
\end{eqnarray*}
where for the last equality we have used the preceeding lemma   with $\lambda=(\mu\mu'^{-1})$ and $ \frac{(\alpha-\delta)}{\alpha}=u\varpi^i.$

From
\begin{eqnarray*}
&&f_{x,z}D_{i,\alpha,\delta}f_{x,z}^{-1}=\\
&&=\begin{pmatrix}\delta+(\delta-\alpha)\frac{\varpi z (x+\varpi z)}{\xi(z)}& \frac{(1+x\varpi z)(x+\varpi z)(\alpha-\delta)}{\xi(z)}\\
\frac{\varpi z(\alpha-\delta)}{\xi(z)}& \alpha+(\alpha-\delta)\frac{\varpi z(x+\varpi z)}{\xi(z)}
\end{pmatrix},
\end{eqnarray*}
only the matrix with $v((\alpha-\delta)\varpi z)\geq l$  and $v((\alpha-\delta)(x+\varpi z)\geq l$ can contribute to $S_f(D_{i,\alpha,\delta}).$ This last condition is also equivalent to $v(\varpi z)\geq l-i$ and $v(x)\geq l-i.$ 
We therefore have:
\begin{eqnarray*}
S_f(D_{i,\alpha,\delta})&=&\sum_{x\in{\mathbb O}_l, z\in {{\mathbb O}_{l-1}}\atop v(x)\geq l-i, v(\varpi z)\geq l-i}
\mu(\delta-(\alpha-\delta)\frac{\varpi z (x+\varpi z)}{\xi(z)})\mu'(\alpha+(\alpha-\delta)\frac{\varpi z(x+\varpi z)}{\xi(z)})\\
=&&  \mu(\delta)\mu'(\alpha)
\sum_{x\in{\mathbb O}_l, z\in {{\mathbb O}_{l-1}} \atop v(x)\geq l-i, v(\varpi z)\geq l-i}\mu(1-\frac{\alpha-\delta}{\delta}\frac{\varpi z (x+\varpi z)}{\xi(z)})\mu'(1+\frac{(\alpha-\delta)}{\alpha}\frac{\varpi z (x+\varpi z)}{\xi(z)})\\
&&=\mu(\delta)\mu'(\alpha)
\sum_{x\in{\mathbb O}_l, z\in {{\mathbb O}_{l-1}} \atop v(x)\geq l-i, v(\varpi z)\geq l-i}\mu(1-\frac{\alpha-\delta}{\delta}\frac{\varpi z (x+\varpi z)}{\xi(z)})\mu'(1+\frac{(\alpha-\delta)}{\delta}\frac{\varpi z (x+\varpi z)}{\xi(z)})\\
&&=\mu(\delta)\mu'(\alpha)\sum_{x\in{\mathbb O}_l, z\in {{\mathbb O}_{l-1}} \atop v(x)\geq l-i, v(\varpi z)\geq l-i}(\mu\mu'^{-1})(1-\frac{\alpha-\delta}{\delta}\frac{\varpi z (x+\varpi z)}{\xi(z)})\\
&&=q^i\mu(\delta)\mu'(\alpha).
\end{eqnarray*}
The proof of the last equality follows the same analysis as the preceeding lemma  with minor adjustments.
Indeed we have to distinguish two cases. If $i\geq l$ then $i=i'+l,$
the summation on $x,$ $z$ being fixed, is $0$ unless $v(\varpi^{i'}\varpi z)\geq l$ and in this case it gives 
$ \mu(\delta)\mu'(\alpha) \vert {\mathbb O}_l \vert (\mu\mu'^{-1})(1-\frac{\alpha-\delta}{\delta}(\frac{(\varpi z)^2}{\xi(z)}).$ From the condition on $i$ we have $(\alpha-\delta)(\varpi z)^2=0,$ therefore $S_f(D_{i,\alpha,\delta})= \mu(\delta)\mu'(\alpha) \vert {\mathbb O}_l \vert\vert \{ z\in {{\mathbb O}_{l-1}}, v(\varpi^{i'}\varpi z)\geq l\}\vert =
\mu(\delta)\mu'(\alpha) \vert {\mathbb O}_l\vert  \vert \{ z\in {\mathbb O}_{l-1}, \varpi z\in \varpi^{l-i'}{\mathbb O}_{i'}\vert =\mu(\delta)\mu'(\alpha) q^i. $
If $i<l$, we proceed as in the previous lemma  setting  $x=\varpi^{l-i}x'$, the summation on $x'$ forces $\varpi z$ to be zero, in order to get a non zero sum.  
As a result  $S_f(D_{i,\alpha,\delta})=\mu(\delta)\mu'(\alpha)\vert {\mathbb O}_i \vert =q^i\mu(\delta)\mu'(\alpha).$

\endproof
$\bullet$ Conjugacy classes of type $B$

\begin{proposition}
$tr(B_{i,\alpha,\beta}\vert_{\Pi_{\mu,\mu'}})=\delta_{i,r-1}q^{r-1}\mu(\alpha)\mu'(\alpha).$
\end{proposition}

\proof
\begin{equation*}
e_{x,y}B_{i,\alpha,\delta}e_{x,y}^{-1}=
\begin{pmatrix}
\alpha+\varpi^i(x-y\varpi \beta(1+xy))& \varpi^{i+1}\beta(1+xy)^2-x^2\varpi^i\\
\varpi^i(1-\beta\varpi y^2)&\alpha+\varpi^{i+1}\beta y(1+xy)-\varpi^i x
\end{pmatrix}.
\end{equation*}

This matrix is in $SK_{l}$ only when  $i\geq l$.   

\begin{equation*}
f_{x,z}B_{i,\alpha,\delta}f_{x,z}^{-1}=
\begin{pmatrix}
\alpha+\frac{\varpi^{i+1}}{\xi(z)}(\beta(\varpi z+x)-z(1+\varpi z x))& \frac{\varpi^{i}}{\xi(z)}((1+\varpi z x)^2-(x+\varpi z)^2\varpi \beta))\\
\frac{\varpi^{i+1}}{\xi(z)}(\beta-\varpi x^2)&\alpha+\frac{\varpi^{i+1}}{\xi(z)}(z(x\varpi z+1)-\beta(x+\varpi z)))
\end{pmatrix}.
\end{equation*}
This matrix is in $SK_{l}$ only when  $i\geq l.$

Therefore if $i<l$ then $S_e(D_{i,\alpha,\delta})=S_f(D_{i,\alpha,\delta})=0.$

If $i\geq l$ we have 
\begin{eqnarray*}
S_e(B_{i,\alpha,\beta})&=&\sum_{x,y\in {\mathbb O_{l}}}
\mu(\alpha+\varpi^i(x-y\varpi \beta(1+xy) ))
\mu'(\alpha+\varpi^{i+1}\beta y(1+xy)-\varpi^i x )\\
&=&\mu(\alpha)\mu'(\alpha)\sum_{x,y\in {\mathbb O_{l}}}
(\mu\mu'{}^{-1})(1+\frac{\varpi^i}{\alpha}(x(1-\varpi y^2\beta)-y\varpi\beta))),
\end{eqnarray*}
the sum over $x$ gives $0$ for $i\leq r-1.$ 

\begin{eqnarray*}
&&S_f(B_{i,\alpha,\beta})=\\
&&\sum_{x\in {\mathbb O_{l}},z\in {\mathbb O}_{l-1}}
\mu(\alpha+\frac{\varpi^{i+1}}{\xi(z)}(\beta(\varpi z+x)-z(1+\varpi z x)))
\mu'(\alpha- \frac{\varpi^{i+1}}{\xi(z)}(\beta(\varpi z+x)-z(1+\varpi z x)))\\
&=&\mu(\alpha)\mu'(\alpha)\sum_{x\in {\mathbb O_{l}},z\in {\mathbb O}_{l-1}}
(\mu\mu'{}^{-1})(1+\frac{\varpi^{i+1}}{\alpha\xi(z)}(\beta(\varpi z+x)-z(1+\varpi z x))).
\end{eqnarray*}
The sum over $x$ gives $0$ unless $i=r-1$  where the result is $\mu(\alpha)\mu'(\alpha)\vert{\mathbb O_{l}}\vert\vert{\mathbb O_{l-1}}\vert=q^{r-1}\mu(\alpha)\mu'(\alpha).$
\endproof

The following table give the complete list of evaluation of characters of principal split representations:

\medskip

\begin{tabular}{|c|c|c|c|c|}\hline
  &$I_{\alpha}$   & $D_{i,\alpha,\delta}$& $C_{i,\alpha,\beta}$ &$B_{i,\alpha,\beta}$ \\ 
\hline\hline
 $ tr({\Pi_{\mu,\mu'}})(.)$ & $q^{r-1}(q+1)\mu(\alpha)\mu'(\alpha)$ &$q^{i}(\mu(\alpha)\mu'(\delta)+ \mu'(\alpha) \mu(\delta) )$ &0&$\delta_{i,r-1}q^{r-1}\mu(\alpha)\mu'(\alpha)$ \\
\hline
\end{tabular}
\bigskip

\subsection{ Cuspidal representations}
Having chosen $\beta$ of type $C_3$, and a lift $\hat{\beta}$  we obtain:
\begin{proposition}
A section of the right cosets of $SK_{l}$  are given by the following set of matrices
 $Y=\{h_{c,d} ,c\in{\mathbb O}_l, d\in{\mathbb   O}_{l}^\times\}$ where 
$ h_{c,d}=\begin{pmatrix}d&0\\c&1\end{pmatrix} .$
\end{proposition}
\proof 
We have $SK_{l}=\{\begin{pmatrix} a+\varpi^l x&b+ \varpi^l y\\-b\hat{\tau}\sigma(\hat{\tau})+\varpi^l z& a+b(\hat{\tau}+\sigma(\hat{\tau}))+\varpi^l t\end{pmatrix}, a,b,x,y, z,t\in { \mathcal O}_r
\}.$ It is an easy exercise to show that the orbits associated to the elements $Y$ are disjoints. Moreover $\vert Y\vert=q^{l}(q^l-q^{l-1})=\vert G/SK_l\vert$ therefore $Y$ is a section of the right cosets of  $SK_{l}$.

\endproof
Frobenius formula therefore gives:

 $$tr(g\vert_{{\mathcal C}_{\nu,\nu'}})=\sum_{t\in Y}(\nu\boxtimes\nu')^0(tgt^{-1})=S_{h}(g,\nu,\nu').$$

$\bullet$ Conjugacy class of type $I.$

$I_{\alpha}$ being central we have

 $tr(I_\alpha\vert_{{\mathcal C}_{\nu,\nu'}})=\frac{\vert G\vert}{\vert SK_{l}\vert}\nu(\alpha)\nu'(\alpha)=q^{r-1}(q-1)\nu(\alpha)\nu'(\alpha).$

$\bullet $
Conjugacy class of type $D.$

\begin{proposition}
$tr(D_{i,\alpha,\delta}\vert_{{\mathcal C}_{\nu,\nu'}})=0$
\end{proposition}

\proof
\begin{equation*}
h_{c,d}D_{i,\alpha,\delta}h_{c,d}^{-1}=
\begin{pmatrix}
\alpha&0\\
(\alpha-\delta)c/d&\delta
\end{pmatrix},
\end{equation*}
therefore 
\begin{eqnarray*}
&&S_{h}(D_{i,\alpha,\delta},\nu,\nu')=
\sum_{c\in{\mathbb O}_l, d\in{\mathbb O}_l^\times
}(\nu\boxtimes\nu')^0(\begin{pmatrix}
\alpha&0\\
(\alpha-\delta)c/d&\delta
\end{pmatrix})\\
&&=\vert {\mathbb O}_l^\times\vert \sum_{c\in{\mathbb O}_l}({\nu\boxtimes\nu'})^0(\begin{pmatrix}
\alpha&0\\
(\alpha-\delta)c&\delta
\end{pmatrix})
\end{eqnarray*}
We have $\alpha-\delta=\varpi^i u$ with $u$ invertible  therefore $h_{c,d}D_{i,\alpha,\delta}h_{c,d}^{-1}\in SK_{l}$ only if $i\geq l$ (this comes from the fact it is lower triangular).

In the case where $i\geq l$ we have
$$\begin{pmatrix}
\alpha&0\\
(\alpha-\delta)c&\delta
\end{pmatrix}=I_{\alpha}\begin{pmatrix}
1&0\\
\frac{\varpi^iuc}{\alpha}&1-\frac{\varpi^iu}{\alpha}
\end{pmatrix}
,$$ therefore

 $$(\nu\boxtimes\nu')(\begin{pmatrix}
\alpha&0\\
(\alpha-\delta)c&\delta
\end{pmatrix})=
\nu(\alpha)\nu'(\alpha)\psi(\frac{1}{\alpha}(\varpi^i uc-Tr(\tau)\varpi^i u)).$$
As a result
\begin{eqnarray*}
S_{h}(D_{i,\alpha,\delta},\nu,\nu')=
\vert {\mathbb O}_l^{\times} \vert \nu(\alpha)\nu'(\alpha)
\sum_{c\in{\mathbb O}_l}\psi(\frac{1}{\alpha}(\varpi^i uc-Tr(\tau)\varpi^i u))=0,
\end{eqnarray*}the summation on $c$ giving  $0$ because $\psi$ is primitive and $i<r.$

\endproof

$\bullet$ Conjugacy class of type ${ B}.$
\begin{proposition}
$tr(B_{i,\alpha,\beta}\vert_{{\mathcal C}_{\nu,\nu'}})=-\delta_{i,r-1}\nu(\alpha)\nu'(\alpha)q^{r-1}.$

\end{proposition}

\proof
We prefer to work with the representative $wB_{i,\alpha,\beta}w^{-1}$ where 
$w=\begin{pmatrix}1&0\\0&1\end{pmatrix},$ it will be easier to compare with results of \cite{LCW}.
We have 
\begin{equation*}
h_{c,d}wB_{i,\alpha,\beta}w^{-1}h_{c,d}^{-1}=
\begin{pmatrix}
\alpha-c\varpi^i&d\varpi^i\\ \frac{\varpi^i}{d}(\varpi\beta-c^2)&\alpha+\varpi^i c
\end{pmatrix}.
\end{equation*}
This matrix belongs to $SK_l$ only when $i\geq l.$

Indeed, 
 if $\begin{pmatrix}
\alpha-c\varpi^i&d\varpi^i\\ \frac{\varpi^i}{d}(\varpi\beta-c^2)&\alpha+\varpi^i c
\end{pmatrix}=\begin{pmatrix}
a+\varpi^l x&b+\varpi^l y \\-bN(\hat{\tau})+\varpi^l y&a+bTr(\hat{\tau})+\varpi^l t
\end{pmatrix}$
we necessarily have $b=d\varpi^i \,\text{mod}\, \varpi^l$, $a=\alpha-\frac{d}{2}Tr(\hat{\tau})\varpi^i \,\text{mod} \,\varpi^l$, 
$c\varpi^i=\frac{d}{2}Tr(\hat{\tau})\varpi^i \text{mod} \varpi^l, \varpi^{i+1}\beta-c^2\varpi^i=-d^2\varpi^iN(\hat{\tau})
\,\text{mod} \,\varpi^l.$
Assuming $i<l$ the last equation implies that $(c^2-d^2N(\hat{\tau}))\varpi^i=0 \,\text{mod}\,\varpi^{i+1},$ but using 
$c=\frac{d}{2}Tr(\hat{\tau} )\text{mod} \varpi^{l-i},$ we necessarily have that $d^{2}((\frac{1}{2}Tr(\hat{\tau} ))^2-N(\hat{\tau}))=d^{2}\tilde{\epsilon}=0 \,\text{mod} \,\varpi$ which contradicts $d$ and $\tilde{\epsilon}$ invertible.

Let us therefore assume that $i\geq l,$ we  have $h_{c,d}wB_{i,\alpha,\beta}w^{-1}h_{c,d}^{-1}\in SK_{l},$ and:

\begin{eqnarray*}
&&S_{h}(B_{i,\alpha,\beta},\nu,\nu')=S_{h}(wB_{i,\alpha,\beta}w^{-1},\nu,\nu')=
\\&&=
\nu(\alpha)\nu'(\alpha)\sum_{c\in{\mathbb O}_l ,d\in{\mathbb O}_l^{\times }  }
\psi(\frac{\varpi^i}{d\alpha}(\varpi\beta-c^2)-\frac{\varpi^i d}{\alpha}N(\tau)+\frac{\varpi^i}{\alpha}Tr(\tau)c)\\
&&= \nu(\alpha)\nu'(\alpha)\sum_{c\in{\mathbb O}_l ,d\in{\mathbb O}_l^{\times } }
\psi(\frac{\varpi^i}{d\alpha}(\varpi\beta-(c-dTr(\tau))^2)+\frac{\varpi^i d}{\alpha}(\frac{\tau-\sigma(\tau))}{2})^2)\\
&&=\nu(\alpha)\nu'(\alpha)\sum_{c\in{\mathbb O}_l ,d\in{\mathbb O}_l^{\times } }
\psi(\frac{\varpi^i}{d\alpha}(\varpi\beta-c^2)+\frac{\varpi^i d\tilde{\epsilon}}{\alpha}) \; (\text{with}\;\tilde{\epsilon}=(\frac{\tau-\sigma(\tau))}{2})^2)\\
&&=\nu(\alpha)\nu'(\alpha)q^{2(i-l)}\sum_{c\in {\mathbb O}_{r-i}, d\in {\mathbb O}_{r-i}^{\times}}
\psi(\frac{\varpi^i}{d\alpha}(\varpi\beta-c^2)+\frac{\varpi^i d\tilde{\epsilon}}{\alpha}).
\end{eqnarray*}
Using Proposition.\ref{evaluationkloosterman1} of the Appendix, because $\varpi\beta$ is not invertible, this sum is equal to $0$ unless $\varpi^{i+1}\beta=0$, i.e $r-i=1$, and the result is equal to $-\nu(\alpha)\nu'(\alpha)q^{r-1}.$

\medskip

$\bullet$ Conjugacy class of type ${ C}.$

Recall that we have denoted  $\Phi\in{\mathcal O}_r^E$ a solution of $\Phi^2={\epsilon}.$
\begin{proposition}
$$tr(C_{i,\alpha,\beta}\vert_{{\mathcal C}_{\nu,\nu'}})=(-q)^i
(\nu(\alpha+\varpi^i\Phi\beta)\nu'(\alpha-\varpi^i\Phi\beta)+\nu(\alpha-\varpi^i\Phi\beta)\nu'(\alpha+\varpi^i\Phi\beta))
$$
\end{proposition}
\proof

We have \begin{equation*}
h_{c,d}C_{i,\alpha,\beta}h_{c,d}^{-1}=
\begin{pmatrix}
\alpha-c\varpi^i\epsilon \beta&d\varpi^i\epsilon\beta\\ \frac{\beta\varpi^i}{d}(1-c^2\epsilon)&\alpha+\varpi^i c\epsilon\beta
\end{pmatrix}.
\end{equation*}
We have to distinguish two cases $i\geq l$ and $i<l.$

$\Diamond$ If $i\geq l$ then this last matrix belongs to $SK_l$ and is equal to  $I_{\alpha}k$ where $k\in K_l.$
Therefore  
\begin{eqnarray*}
&&S_{h}(C_{i,\alpha,\beta},\nu,\nu')=
\nu(\alpha)\nu'(\alpha)\sum_{c\in{\mathbb O}_l, d\in {\mathbb O}_l^\times}
\psi(\frac{\beta\varpi^i}{d\alpha}(1-c^2\epsilon)-\frac{\varpi^i d\epsilon\beta}{\alpha}N(\tau)+
\frac{\varpi^i}{\alpha}Tr(\tau)c\epsilon\beta)\\
&&= \nu(\alpha)\nu'(\alpha)\sum_{c\in{\mathbb O}_l,d\in {\mathbb O}_l^\times }
\psi(\frac{\beta\varpi^i}{d\alpha}(1-(c+\frac{d}{2}Tr(\tau))^2\epsilon)-\frac{\varpi^i d\epsilon\beta}{\alpha}N(\tau)+\frac{\varpi^i}{\alpha}Tr(\tau)(c+\frac{d}{2}Tr(\tau))\epsilon\beta)\\
&&=\nu(\alpha)\nu'(\alpha)\sum_{c\in{\mathbb O}_l,d\in {\mathbb O}_l^\times }
\psi(\frac{\varpi^i\beta}{d\alpha}(1-c^2\epsilon)+\frac{\varpi^i \epsilon \beta d\tilde{\epsilon}}{\alpha})\\
&&=\nu(\alpha)\nu'(\alpha)
\sum_{c\in{\mathbb O}_l,d\in {\mathbb O}_l^\times }
\psi(\frac{\varpi^i\beta}{d\alpha}(\epsilon-c^2\epsilon^2)+\frac{\varpi^i \beta d\tilde{\epsilon}}{\alpha})\\
&&=\nu(\alpha)\nu'(\alpha)
\sum_{c\in{\mathbb O}_l,d\in {\mathbb O}_l^\times }
\psi(\frac{\varpi^i\beta}{d\alpha}(\epsilon-c^2)+\frac{\varpi^i \beta d\tilde{\epsilon}}{\alpha})\\
&&=\nu(\alpha)\nu'(\alpha)q^{2(i-l)}\sum_{c\in {\mathbb O}_{r-i}, d\in {\mathbb O}_{r-i}^{\times}}
\psi(\frac{\varpi^i\beta}{d\alpha}(\epsilon-c^2)+\frac{\varpi^i \beta d\tilde{\epsilon}}{\alpha}).
\end{eqnarray*}
Applying  Proposition (\ref{evaluationkloosterman1}) of the Appendix (evaluation of Sali\'e sums), this is equal to 
\begin{eqnarray*}
tr(C_{i,\alpha,\beta}\vert_{{\mathbb C}_{\nu,\nu'}})&=&\nu(\alpha)\nu'(\alpha)q^{2(i-l)}(-q)^{r-i}
(\psi(\frac{\varpi^i\beta}{\alpha}2u)+\psi(-\frac{\varpi^i\beta}{\alpha}2u))\\
&&(\text{where} \;u^2=\epsilon\tilde{\epsilon}, 2u=\Phi(\tau-\sigma(\tau)))\\
&&=\nu(\alpha)\nu'(\alpha)(-q)^i
(\psi(\frac{\varpi^i\beta}{\alpha}\Phi(\tau-\sigma(\tau)))+\psi(-\frac{\varpi^i\beta}{\alpha}\Phi(\tau-\sigma(\tau))),\\
&&=\nu(\alpha)\nu'(\alpha)(-q)^i(\nu(1+\varpi^i\Phi\frac{\beta}{\alpha})\nu'(1-\varpi^i\Phi\frac{\beta}{\alpha})+\nu(1-\varpi^i\Phi\frac{\beta}{\alpha})\nu'(1+\varpi^i\Phi\frac{\beta}{\alpha}))\\
&&=(-q)^i
(\nu(\alpha+\varpi^i\Phi\beta)\nu'(\alpha-\varpi^i\Phi\beta)+\nu(\alpha-\varpi^i\Phi\beta)\nu'(\alpha+\varpi^i\Phi\beta))
\end{eqnarray*}
which gives the announced result. 

Remark: This corrects misprints in \cite{LCW} (Page 1306 it should be $p^{i}$ and not $p^{l-i}$ in the expression of $tr(C_{i,\alpha,\beta}\vert_{{\mathcal C}_{\nu,\nu'}})$ and 
moreover they have  considered the case were $\epsilon=\tilde{\epsilon}=u)$.

$\Diamond$ We now proceed with $i<l.$ The same method as developped in \cite{LCW} can be applied and this is a non trivial result.
We can ask what are the $(c,d)\in {\mathbb O}_l\times  {\mathbb O}_l^\times$ such that $h_{c,d}C_{i,\alpha,\beta}h_{c,d}^{-1}\in SK_l$, this condition is equivalent to the existence 
of  $X,Y,Z,T\in {\mathcal O}_r $ such that 
$$h_{c,d}C_{i,\alpha,\beta}h_{c,d}^{-1}=\begin{pmatrix} X&Y\\- Y\hat{\tau} \sigma(\hat{\tau})+\varpi^l Z&X+Y(\hat{\tau}+\sigma(\hat{\tau}))+\varpi^l T
 \end{pmatrix}.$$
This imply $Y=d\varpi^i\epsilon\beta$, $2c\varpi^i\epsilon\beta=Y Tr(\hat{\tau}) +\varpi^l T.$ Therefore $c=\frac{d}{2}Tr(\hat{\tau})+\varpi^{l-i}T, T\in{\mathbb O}_i.$
From $\frac{\beta\varpi^i}{d}(1-c^2\epsilon)=-YN(\tau)+\varpi^l Z$ , a little algebra implies that $d=\pm u^{-1}+d' \varpi^{l-i}, d'\in {\mathbb O}_i,$ and 
$c=\pm \frac{u^{-1}}{2}Tr(\hat{\tau})+c'\varpi^{l-i}, c'\in {\mathbb O}_i.$
Inversely this is a sufficient condition  on $c,d$ for having $h_{c,d}C_{i,\alpha,\beta}h_{c,d}^{-1}\in SK_l.$

We now fix $d= u^{-1}+d' \varpi^{l-i}, c=\frac{u^{-1}}{2}Tr(\tau)+c'\varpi^{l-i}$, the other choice of sign follows the same method.

Under this condition, we now factorise: 
$h_{c,d}C_{i,\alpha,\beta}h_{c,d}^{-1}=sk, s\in S, k\in K_l.$ This is not unique, we can choose $s=\begin{pmatrix} a&b\\-b\hat{\tau} \sigma(\hat{\tau})&a+b(\hat{\tau}+\sigma(\hat{\tau}))
 \end{pmatrix}$, with $a=\alpha-\frac{u^{-1}}{2}Tr(\hat{\tau})\varpi^i\epsilon\beta, b=u^{-1}\varpi^i\epsilon\beta.$
The matrix $k=s^{-1}h_{c,d}C_{i,\alpha,\beta}h_{c,d}^{-1}=\begin{pmatrix}1+\varpi^l x&\varpi^l y\\ \varpi^l z&1+\varpi^l t
\end{pmatrix}.$
Note that $det(s)=a^2+abTr(\hat{\tau})+b^2N(\hat{\tau})$ which after simplifications gives $det(s)=\alpha^2-\varpi^{2i}\epsilon\beta^2.$
Therefore the matrix $k$ can be computed as 
$$k=\frac{1}{\alpha^2-\varpi^{2i}\epsilon\beta^2}
\begin{pmatrix} a&b\\-b\hat{\tau} \sigma(\hat{\tau})&a+b(\hat{\tau}+\sigma(\hat{\tau}))
 \end{pmatrix}\begin{pmatrix}
\alpha-c\varpi^i\epsilon \beta&d\varpi^i\epsilon\beta\\ \frac{\beta\varpi^i}{d}(1-c^2\epsilon)&\alpha+\varpi^i c\epsilon\beta
\end{pmatrix}.$$

From this equation we obtain, after a direct lengthy computation 
\begin{eqnarray*}
&&(\alpha^2-\varpi^{2i}\epsilon\beta^2)(\varpi^l z-\varpi^l y N(\tau)+Tr(\tau )(1+\varpi^l t))=\\
&&=2b\alpha N(\tau)+\frac{a\beta\varpi^i}{d}(1-c^2\epsilon)-ad\varpi^i \epsilon\beta N(\tau)+a(\alpha+\varpi^i c\epsilon\beta)Tr(\tau).
\end{eqnarray*}
This expression can be simplified by different changes of variables.
Defining $\tilde{c}=c-\frac{d}{2}Tr(\tau)$, we obtain 
\begin{eqnarray*}
&&d(\alpha^2-\varpi^{2i}\epsilon\beta^2)(\varpi^l z-\varpi^l y N(\tau)+Tr(\tau )\varpi^l t)=\\
&&=2bd\alpha N(\tau)+ad\alpha Tr(\tau )+a\beta\varpi^i (1-\tilde{c}^2\epsilon)+ad^2\varpi^i \beta \epsilon\tilde{\epsilon }-d(\alpha^2-\varpi^{2i}\epsilon\beta^2)Tr(\tau)\\
&&=a d^2\varpi^i u^2\beta-2u\varpi^i\beta a d+\alpha\beta\varpi^i-\tilde{c}^2\epsilon a\beta\varpi^i.
\end{eqnarray*}
Setting $d=u^{-1}+f\varpi^{l-i},\tilde{c}=e\varpi^{l-i}$ with $e,f\in{\mathbb O }_i$ in the last expression, we obtain:
$$\varpi^l z-\varpi^l y N(\tau)+Tr(\tau )\varpi^l t=
\frac{1}{\alpha^2-\varpi^{2i}\epsilon\beta^2}\frac{ua\beta\epsilon\pi^{2l-i}}{1+fu\varpi^{l-i}}(f^2\tilde{\epsilon}-e^2).
$$
As a result we obtain for  $d=u^{-1}+f\varpi^{l-i},\tilde{c}=e\varpi^{l-i}:$
\begin{equation*}
(\nu\boxtimes\nu')(h_{c,d}C_{i,\alpha,\beta}h_{c,d}^{-1})=\nu(a+b\hat{\tau})\nu'(a+b\sigma(\hat{\tau}))
\psi(\varpi^l z-\varpi^l y N(\tau)+Tr(\tau )\varpi^l t)
\end{equation*}
We can simplify:
\begin{eqnarray*}
&&\nu(a+b\hat{\tau})\nu'(a+b\sigma(\hat{\tau}))=\\
&&=\nu(\alpha-\frac{u^{-1}}{2}\varpi^i\epsilon\beta Tr(\hat{\tau})+u^{-1}\varpi^i\epsilon\beta\hat{\tau})
\nu'(\alpha-\frac{u^{-1}}{2}\varpi^i\epsilon\beta Tr(\hat{\tau})+u^{-1}\varpi^i\epsilon\beta\sigma(\hat{\tau}))\\
&&=\nu(\alpha+\frac{\varpi^i\epsilon(\tau-\sigma(\tau))\beta}{2u} )\nu'(\alpha-\frac{\varpi^i\epsilon(\tau-\sigma(\tau))\beta}{2u})\\
&&=\nu(\alpha+\varpi^i\Phi\beta)\nu'(\alpha-\varpi^i\Phi\beta).
\end{eqnarray*}
As a result 
\begin{equation*}
S_{h}(C_{i,\alpha,\beta},\nu,\nu')=\nu(\alpha+\varpi^i\Phi\beta)\nu'(\alpha-\varpi^i\Phi\beta)S^+  +\nu(\alpha-\varpi^i\Phi\beta)\nu'(\alpha+\varpi^i\Phi\beta)S^-\\
\end{equation*}
where 
 $$S^{\pm}=\sum_{e,f\in {\mathbb O}_i}{}\psi(\pm\frac{1}{\alpha^2-\varpi^{2i}\epsilon\beta^2}\frac{ ua\beta\epsilon\pi^{2l-i}}{1\pm fu\varpi^{l-i}}(f^2\tilde{\epsilon}-e^2)).$$
Let us define  the characters  $\lambda^{\pm}:{\mathcal O}_i\rightarrow{\mathbb C}^\times,\lambda^\pm(z)=\psi(\frac{\pm ua\beta\epsilon\pi^{2l-i}z}{\alpha^2-\varpi^{2i}\epsilon\beta^2}).$

We have 
\begin{eqnarray*}
&&S^\pm=\sum_{e,f\in {\mathcal O}_i}\lambda^\pm(\frac{f^2\tilde{\epsilon}-e^2)}{1\pm fu\varpi^{l-i}})\\
&&=\sum_{e',f'\in {\mathcal O}_i}\lambda^\pm(f'^2\tilde{\epsilon}-e'^2)\;\text{with } \;e'=\frac{e}{\sqrt{1\pm fu\varpi^{l-i}}},  f'=\frac{f}{\sqrt{1\pm fu\varpi^{l-i}}}\\
&&=(-q)^i.
\end{eqnarray*}
We have used to conclude the proposition \ref{KloostermanGauss} of the appendix using properties of Gauss sums.
\endproof
As a result 
\begin{eqnarray*}
S_{h}(C_{i,\alpha,\beta},\nu,\nu')=(-q)^i(\nu(\alpha+\varpi^i\Phi\beta)\nu'(\alpha-\varpi^i\Phi\beta)+\nu(\alpha-\varpi^i\Phi\beta)\nu'(\alpha+\varpi^i\Phi\beta)),
\end{eqnarray*}
which is the desired result.
\endproof

\medskip
The following table give the complete list of evaluation of characters of cuspidal  representations:
\medskip

\begin{tabular}{|c|c|c|c|c|}\hline
  &$I_{\alpha}$   & $D_{i,\alpha,\delta}$& $C_{i,\alpha,\beta}$ &$B_{i,\alpha,\beta}$ \\ 
\hline\hline
$tr(.\vert_{{\mathcal C}_{\nu,\nu'}})$ & $q^{r-1}(q-1)\nu(\alpha)\nu'(\alpha)$
 &$0$ &$\begin{array}{c}(-q)^i
(\nu(\alpha+\varpi^i\Phi\beta)\nu'(\alpha-\varpi^i\Phi\beta)\\
+\nu(\alpha-\varpi^i\Phi\beta)\nu'(\alpha+\varpi^i\Phi\beta))
\end{array}$&$-\delta_{i,r-1}\nu(\alpha)\nu'(\alpha)q^{r-1}$ \\
\hline
\end{tabular}
\bigskip

\subsection{ Non principal split representations}

Let $\Delta\in {\mathbb O}_l\cap {\mathfrak M},$ let $\hat{\Delta}$ a lift in ${\mathcal O}_r$ of $\Delta.$
We have
 $T(\psi_{\beta(C_4'(\Delta,0))})=SK_{l}$ with $S=S(\hat{\Delta}, 0)=
\{\begin{pmatrix}a&b\\
-\hat{\Delta} b& a\end{pmatrix}, a\in  {\mathcal O}_r^\times,  b\in {\mathcal O}_r\}.$

\begin{proposition}
A section of the right cosets of $SK_{l}$  is  given by the following set of matrices $Y\cup Z$ where 
 $Y=\{h_{c,d} ,c\in{\mathbb O}_l, d\in{\mathbb  O}_{l}^\times\}$ with
$ h_{c,d}=\begin{pmatrix}d&0\\c&1\end{pmatrix} $ and $Z=\{ h_{\varpi c,d}w,  c\in{\mathbb O}_{l-1}, d\in{\mathbb  O}_{l}^\times\}$ and $w=\begin{pmatrix}0&1\\1&0\end{pmatrix}.$
\end{proposition}
\proof

We use  $SK_{l}=\{\begin{pmatrix}X&Y\\-\hat{\Delta}Y+\varpi^l Z&X+\varpi^l T\end{pmatrix}, X\in {\mathcal O}_r^{\times}, Y,Z,T\in {\mathcal O}_r\}$.  It is a direct verification.
\endproof
\medskip

An important remark: note that $w$ is put on the right  of $h_{\varpi c,d}$ and this is essential for the proposition to be true. 
In \cite{LCW} an analog of this proposition is stated in Lemma 6.1  and in  section 6.2  page 1315 but $w$ is incorrectly put on the left and the corresponding  set of
 elements $E_{cd}$ and $F_{cd}$ (in their notations) do not provide a set of representative of right cosets.
 Indeed  with the notations of Lemma 6.1, if we define  the group $L=\{\begin{pmatrix}x& \varpi^j y\beta\\y&x \end{pmatrix}\},$  we have 
$F_{cd}F_{c'd'}^{-1}=\begin{pmatrix}dd'{}^{-1}&0\\\frac{(c-c')\varpi}{d'}&1 \end{pmatrix}$, this matrix belongs to $L$ as soon as $d=d', v(c-c')\geq i-j-1$ ($j\leq i)$ with their notations. Therefore there exist  $c, c'$ such that $F_{cd}\not=F_{c'd}$  and  $LF_{cd}=LF_{c'd}$ contradicting the Lemma 6.1.
Our remark applies as well  for the section 6.2.
As a consequence the results obtained in \cite{LCW} concerning the value of the characters of the split non principal representations  cannot be trusted. 
\medskip

$\bullet$ Conjugacy class of type $I.$

$I_{\alpha}$ being central we have

 $tr(I_\alpha\vert_{{\Xi}_{a,\Delta,\theta}})=\frac{\vert G\vert}{\vert SK_{l}\vert} \tilde{\chi}_a(\alpha^2)\theta(I_\alpha)=(q^2-1)q^{r-2}\tilde{\chi}_a(\alpha^2)\sigma(\alpha),$
where $\sigma$ is the multiplicative character associated to $\theta$.

$\bullet$ Conjugacy class of type $D.$

We have $h_{c,d}D_{i,\alpha,\delta}h_{c,d}^{-1}=\begin{pmatrix}\alpha&0\\c(\alpha-\delta)/d&\delta\end{pmatrix}$ and 
$h_{\varpi c,d}wD_{i,\alpha,\delta}w h_{\varpi c,d}^{-1}=\begin{pmatrix}\delta&0\\\varpi c(\delta-\alpha)/d&\alpha\end{pmatrix},$ as a result the conjugacy class of $D_{i,\alpha,\delta}$ intersects $ SK_{l}$ only  when $i\geq l$  because the valuation of the difference of the matrix elements on the diagonal has to be greater or equal to $l.$

Therefore we obtain that $tr(D_{i,\alpha,\delta}\vert_{{\Xi}_{a,\Delta,\theta}})=0$ when $i<l.$

When $i\geq l$ we have both $h_{c,d}D_{i,\alpha,\delta}h_{c,d}^{-1}\in SK_{l}, h_{\varpi c,d}wD_{i,\alpha,\delta}w h_{\varpi c,d}^{-1}\in SK_{l}$. Noting that $h_{c,d}D_{i,\alpha,\delta}h_{c,d}^{-1}=
I_{\alpha}\begin{pmatrix}1&0\\c(\alpha-\delta))/(d\alpha)&1-(\alpha-\delta)/\alpha\end{pmatrix}$ and similarly for  $h_{\varpi c,d}wD_{i,\alpha,\delta}w h_{\varpi c,d}^{-1}$ we obtain 
\begin{eqnarray*}
tr(D_{i,\alpha,\delta}\vert_{{\Xi}_{\Delta,\theta}})=
\sum_{c\in {\mathbb O}_l, d\in {\mathbb O}_l^{\times}}\theta(I_{\alpha})\psi(\frac{c}{d}\frac{\alpha-\delta}{\alpha})+
\sum_{ c\in {\mathbb O}_{l-1}, d\in {\mathbb O}_l^{\times}}\theta(I_{\delta})\psi(\frac{\varpi c}{d}\frac{\delta-\alpha}{\delta}).
\end{eqnarray*}
The first sum is always $0$, because $\alpha-\delta\in \varpi^i{\mathcal O}_r$ with  $i<r,$ and therefore for fixed $d$ the sum over $c$ gives $0.$
The second sum is zero for the same reason unless $i=r-1.$ In the case $i=r-1$ we have $\varpi(\alpha-\delta)=0,$ therefore the second sum
 is equal to $\theta(I_{\delta})\vert {\mathbb O}_l^{\times}\vert \vert  {\mathbb O}_{l-1}\vert=
\theta(I_{\delta})(q-1)q^{r-2}.$ 

We finally  obtain 

$tr(D_{i,\alpha,\delta}\vert_{{\Xi}_{a,\Delta,\theta}})=\delta_{i,r-1}\tilde{\chi}_a(\alpha)\tilde{\chi}_a(\delta)\sigma(\delta)(q-1)q^{r-2}.$ 

Note that because the restriction of $\sigma$ to $1+\varpi^l {\mathcal O}_r$ is trivial we have $ \sigma(\alpha)=\sigma(\delta)$ when $i=r-1,$ the result is symmetric in the exchange of $\alpha$ and $\delta,$ as it should be.
\medskip

$\bullet$ Conjugacy class of type $C.$

We have 
\begin{equation*}h_{c,d}C_{i,\alpha,\beta}h_{c,d}^{-1}=
\begin{pmatrix}\alpha-c\varpi^i\epsilon\beta&d\varpi^i\epsilon\beta\\
\beta\varpi^i(1-c^2\epsilon)/d&\alpha+c\varpi^i\epsilon\beta
\end{pmatrix}
\end{equation*}
this matrix does not belong to $SK_{l}$ when $i<l.$ 

Indeed if $h_{c,d}C_{i,\alpha,\beta}h_{c,d}^{-1}=
\begin{pmatrix}X&Y\\-\hat{\Delta}Y+\varpi^l Z&X+\varpi^l T
\end{pmatrix}$
then $Y=d\varpi^i\epsilon\beta$ is of valuation $i$, but  $-\hat{\Delta}Y+\varpi^l Z$ is of valuation strictly bigger to $i$ because $\hat{\Delta}$ is of positive valuation and $i<l,$ therefore one cannot have $-\hat{\Delta}Y+\varpi^l Z=\beta\varpi^i(1-c^2\epsilon)/d$ which is of valuation $i.$

We have 
\begin{equation*}h_{\varpi c,d}w C_{i,\alpha,\beta} wh_{\varpi c,d}^{-1}=
\begin{pmatrix}\alpha-\varpi c\varpi^i\beta&d\varpi^i\beta\\
\beta\varpi^i(\epsilon-(\varpi c)^2)/d&\alpha+\varpi c\varpi^i\beta
\end{pmatrix},
\end{equation*}
for the same reason if $i<l$ then $h_{\varpi c,d}w C_{i,\alpha,\beta} wh_{\varpi c,d}^{-1}\notin SK_l.$

As a result we get
 $tr(C_{i,\alpha,\beta}\vert_{{\Xi}_{a,\Delta,\theta}})=0$ when $i<l.$

\medskip
If  $i\geq l,$ we have 
\begin{equation*}h_{c,d}C_{i,\alpha,\beta}h_{c,d}^{-1}=I_{\alpha}
\begin{pmatrix}1-\frac{c}{\alpha}\varpi^i\epsilon\beta&\frac{d}{\alpha}\varpi^i\epsilon\beta\\
\frac{\beta}{\alpha d}\varpi^i(1-c^2\epsilon)&1+\frac{c}{\alpha}\varpi^i\epsilon\beta
\end{pmatrix}\in SK_{l}.
\end{equation*}

Therefore we obtain 
\begin{eqnarray*}
&&tr(C_{i,\alpha,\beta}\vert_{{\Xi}_{\Delta,\theta}})=\\
&&=\sum_{c\in {\mathbb O}_l, d\in {\mathbb O}_l^{\times}}\theta(I_{\alpha})\psi(\frac{\beta\varpi^i}{\alpha d}(1-c^2\epsilon)-{\Delta} \frac{d}{\alpha}\varpi^i\epsilon\beta)+
\sum_{ c\in {\mathbb O}_{l-1}, d\in {\mathbb O}_l^{\times}}\theta(I_{\alpha})
\psi(\frac{\beta\varpi^i}{\alpha d}(\epsilon-(\varpi c)^2)-{\Delta} \frac{d}{\alpha}\varpi^i\beta)\\
&&=q^{2(i-l)}\sum_{c\in {\mathbb O}_{r-i}, d\in {\mathbb O}_{r-i}^{\times}}\theta(I_{\alpha})\psi(\frac{\beta\varpi^i}{\alpha d}(1-c^2\epsilon)-{\Delta} \frac{d}{\alpha}\varpi^i\epsilon\beta)+\\
&&+q^{2(i-l)}\sum_{ c\in {\mathbb O}_{r-i-1}, d\in {\mathbb O}_{r-i}^{\times}}\theta(I_{\alpha})
\psi(\frac{\beta\varpi^i}{\alpha d}(\epsilon-(\varpi c)^2)-{\Delta} \frac{d}{\alpha}\varpi^i\beta).
\end{eqnarray*}
In order to compute the first term of this  sum we remark that $1-c^2\epsilon$ is invertible and that the map $d\mapsto \frac{1}{d}-d{\Delta}u$ from ${\mathcal O}_{r-i}^{\times}$  to ${\mathcal O}_{r-i}^{\times},$ where $u\in {\mathcal O}_{r-i},$ is bijective. 
As a result we have 
\begin{eqnarray*}
&&\sum_{c\in {\mathbb O}_{r-i}, d\in {\mathbb O}_{r-i}^{\times}}\psi(\frac{\beta\varpi^i}{\alpha d}(1-c^2\epsilon)-{\Delta}
 \frac{d}{\alpha}\varpi^i\epsilon\beta)=\\
&&=\sum_{c\in {\mathbb O}_{r-i}, d\in {\mathbb O}_{r-i}^{\times}}\psi(\frac{\beta\varpi^i}{\alpha d}(1-c^2\epsilon))\\
&&=\sum_{c\in {\mathbb O}_{r-i}, d\in {\mathbb O}_{r-i}^{\times}}\psi(\frac{\beta\varpi^i}{\alpha}(1-c^2\epsilon)d)\\
&&=\sum_{c\in {\mathbb O}_{r-i}, d\in {\mathbb O}_{r-i}}\psi(\frac{\beta\varpi^i}{\alpha}(1-c^2\epsilon)d)-
\sum_{c\in {\mathbb O}_{r-i},  d\in {\mathbb O}_{r-i-1}}\psi(\frac{\beta\varpi^i}{\alpha}(1-c^2\epsilon)\varpi d)
\end{eqnarray*}
Fixing $c$ the sum over $d$ gives $0$ when $r-i>1.$
When $i=r-1$, we have to evaluate $\sum_{c\in {\mathbb O}_{1}, d\in {\mathbb O}_{1}^{\times}}\psi(\frac{\beta\varpi^i}{\alpha}(1-c^2\epsilon)d).$ For fixed $c$ the sum over $d$ gives $-1,$ therefore the value of this sum is equal to $-\vert  {\mathbb O}_{1}\vert=-q.$
We obtain that 
\begin{equation*}
\sum_{c\in {\mathbb O}_l, d\in {\mathbb O}_l^{\times}}\theta(I_{\alpha})\psi(\frac{\beta\varpi^i}{\alpha d}(1-c^2\epsilon)-\hat{\Delta} \frac{d}{\alpha}\varpi^i\epsilon\beta)=-\theta(I_{\alpha})q^{r-1}\delta_{i,r-1}.
\end{equation*}
The second sum is evaluated with the same technique and we obtain
\begin{equation*}
\sum_{ c\in {\mathbb O}_{l-1}, d\in {\mathbb O}_l^{\times}}\theta(I_{\alpha})
\psi(\frac{\beta\varpi^i}{\alpha d}(\epsilon-(\varpi c)^2)-\hat{\Delta} \frac{d}{\alpha}\varpi^i\beta)=-q^{r-2}\theta(I_{\alpha})\delta_{i,r-1}.
\end{equation*}
As a result we get:
\begin{equation*}
tr(C_{i,\alpha,\beta}\vert_{{\Xi}_{\Delta,\theta}})=-q^{r-2}(q+1)\theta(I_{\alpha})\delta_{i,r-1}=-q^{r-2}(q+1)\sigma(\alpha)\delta_{i,r-1}.
\end{equation*}
Finally:
$$tr(C_{i,\alpha,\beta}\vert_{{\Xi}_{a,\Delta,\theta}})=
-\tilde{\chi}_a(det(C_{i,\alpha,\beta}))\sigma(\alpha)q^{r-2}(q+1)\delta_{i,r-1}.$$

Because $\alpha^2-\pi^{2i}\epsilon\beta^2=\alpha^2$ when $i=r-1$, we finally obtain:
$$tr(C_{i,\alpha,\beta}\vert_{{\Xi}_{a,\Delta,\theta}})=
-\tilde{\chi}_a(\alpha^2)\sigma(\alpha)q^{r-2}(q+1)\delta_{i,r-1}.$$

$\bullet$ Conjugacy class of type $B.$

The case of conjugacy class of type $B$ is much more involved than the other conjugacy classes. This comes from two difficulties: there is no neat description of the non principal split representations and moreover  the character depends on  three positive integers $i, j, k.$ We have done a careful analysis of this case but there are cases where we cannot give explicit closed formulas.

Note also that there are numerous mistakes in the analysis of \cite{LCW} concerning this case: the given set of representatives of right cosets is not a set of representatives as already mentionned, the evaluation  of the sum denoted $P$ of \cite{LCW} is mistaken as  noted  and corrected in  \cite{Ma} (but only for the case ${\mathbb Z}_p$) and moreover  they have used the same parameter $\beta$ for the parametrisation of the conjugacy classes and also for the parametrisation of the representation  ${\Xi}_{a,\Delta,\theta}$ (i.e they have called our $\Delta$ also $\beta$). Therefore their  evaluation of the characters, which was nevertheless not given in a closed form for  some cases in $i,j,k$    cannot be trusted for these conjugacy classes.
\medskip

We have 
$$h_{c,d}B_{i,\alpha,\beta}h_{c,d}^{-1}=\begin{pmatrix}\alpha-\varpi^{i+1}\beta&d\varpi^{i+1}\beta\\(\varpi^i-c^2\varpi^{i+1}\beta)/d&\alpha+c\varpi^{i+1}\beta\end{pmatrix},$$ 
$$h_{\varpi c,d}wB_{i,\alpha,\beta}w h_{\varpi c,d}^{-1}=
\begin{pmatrix}\alpha-c\varpi^{i+1}&d\varpi^i\\
(\varpi^{i+1}\beta-\varpi^{i+2}c^2)/d&\alpha+c\varpi^{i+1}\end{pmatrix}.$$

At this point let us define $j=v(\beta)$, we have $0\leq j\leq r-1$ and let $k=v(\hat{\Delta})$ we have $1\leq k\leq r.$ Note that $1+i+j\leq r.$
It will be convenient to chose $\beta'$ and $\hat{\Delta}'$ inversible such that $\beta=\varpi^j\beta', \hat{\Delta}=\varpi^k\hat{\Delta}'.$

\begin{lemma}
$h_{c,d}B_{i,\alpha,\beta}h_{c,d}^{-1}$ belongs to $SK_l$ only if $i\geq l.$ 
\end{lemma}
\proof
Assume that $h_{c,d}B_{i,\alpha,\beta}h_{c,d}^{-1}\in SK_l$, then there exists $Y,Z\in {\mathcal O}_r$  such that $Y=d\varpi^{i+1}\beta$, $(\varpi^i-c^2\varpi^{i+1}\beta)/d=-\hat{\Delta}Y+\varpi^l Z$.  Therefore  $v(Y)=i+j+1$ and $i=v((\varpi^i-c^2\varpi^{i+1}\beta)/d)=v(-\hat{\Delta}Y+\varpi^l Z).$
We have $v(\hat{\Delta}Y)=i+j+k+1$ and $v(\varpi^l Z)=l+v(Z).$
As a result $i\geq min (i+j+k+1,l+v(y))$ which is possible only if $i\geq l.$

\endproof

We consider as usual two cases. 

$\Diamond$ The first case is when $i\geq l,$ in this case both $h_{c,d}B_{i,\alpha,\beta}h_{c,d}^{-1}$ and $h_{\varpi c,d}wB_{i,\alpha,\beta}w h_{\varpi c,d}^{-1}$ belong to $SK_l.$
From $h_{c,d}B_{i,\alpha,\beta}h_{c,d}^{-1}=
I_\alpha\begin{pmatrix}1-\varpi^{i+1}\beta/\alpha&d\varpi^{i+1}\beta/\alpha\\(\varpi^i-c^2\varpi^{i+1}\beta)/(\alpha d)&1+c\varpi^{i+1}\beta/\alpha\end{pmatrix},$ the contribution of the elements $h_{c,d}B_{i,\alpha,\beta}h_{c,d}^{-1}$ to the character of the representation ${\Xi}_{\Delta,\theta} $ is given by the following sum:
\begin{eqnarray*}
&&\chi_1=\sum_{c\in {\mathbb O}_l, d\in {\mathbb O}_l^{\times}}\theta(I_{\alpha})
\psi(\frac{\varpi^i}{\alpha}(\frac{1-c^2\varpi\beta}{ d}-\hat{\Delta}d\varpi\beta))\\
&&=q^{2(i-l)}\theta(I_{\alpha})\sum_{c\in {\mathbb O}_{r-i}, d\in {\mathbb O}_{r-i}^{\times}}
\psi(\frac{\varpi^i}{\alpha}(\frac{1-c^2\varpi\beta}{ d}-\hat{\Delta}d\varpi\beta)).
\end{eqnarray*}
Using the fact that  $d\mapsto \frac{1}{d}-d\hat{\Delta}u$ from ${\mathcal O}_{r-i}^{\times}$  to ${\mathcal O}_{r-i}^{\times},$ where $u\in {\mathcal O}_{r-i}$ is a bijection,  the evaluation of this sum follows the same procedure as in the case of conjugacy classes of type $C$ and we obtain:
$$\chi_1=-\theta(I_\alpha)q^{r-1}\delta_{i,r-1}.$$

From 
$h_{\varpi c,d}wB_{i,\alpha,\beta}w h_{\varpi c,d}^{-1}=I_{\alpha}
\begin{pmatrix}1-c\varpi^{i+1}/\alpha&d\varpi^i/\alpha\\
(\varpi^{i+1}\beta-\varpi^{i+2}c^2)/(\alpha d)&1+c\varpi^{i+1}/\alpha\end{pmatrix}$ the contribution of the elements $h_{\varpi c,d}wB_{i,\alpha,\beta}w h_{\varpi c,d}^{-1}$ to the character of the representation ${\Xi}_{\Delta,\theta} $ is given by the following sum:

\begin{eqnarray*}
&&\chi_2(i,j,k)=\theta(I_{\alpha})\sum_{ c\in {\mathbb O}_{l-1}, d\in {\mathbb O}_l^{\times}}
\psi(\frac{\varpi^{i+1}\beta-\varpi^{i+2}c^2}{\alpha d}-\hat{\Delta}\frac{d}{\alpha}\varpi^i\beta)\\
&&=\theta(I_{\alpha})\sum_{ c\in {\mathbb O}_{l-1}, d\in {\mathbb O}_l^{\times}}
\psi(\frac{\varpi^i}{\alpha}(\frac{\varpi^{j+1}\beta'-\varpi^{2}c^2}{d}-d\varpi^{j+k}\hat{\Delta}'\beta')=\\
&&=q^{2(i-l)}\theta(I_{\alpha})\sum_{ c\in {\mathbb O}_{r-i-1}, d\in {\mathbb O}_{r-i}^{\times}}
\psi(\frac{\varpi^i}{\alpha}(\frac{\varpi^{j+1}\beta'-\varpi^{2}c^2}{d}-d\varpi^{j+k}\hat{\Delta}'\beta').
\end{eqnarray*}
Note that when $i=r-1$, $\psi$ is evaluated on the $0$ element, therefore $\chi_2=\theta(I_\alpha)\vert{\mathbb O}_{l-1}\vert
\vert{\mathbb O}_{l}^{\times}\vert=\theta(I_\alpha)q^{r-2}(q-1).$ Therefore when $i=r-1$, we obtain that the character of  ${\Xi}_{\Delta,\theta} $ is 
$\chi_1+\chi_2=-\theta(I_\alpha)q^{r-2}.$

When $i<r-1$, the value of the character reduces to $\chi_2(i,j,k).$ This sum 
 can be simplified as explained by  \cite{Ma} in the case of $\mathbb {Z}_p$, the  generalisation to  the case of ${\mathcal O}$ is provided in the appendix of our work.

$\Diamond$ The second case is when $i<l.$

We have 
$h_{\varpi c,d}wB_{i,\alpha,\beta}w h_{\varpi c,d}^{-1}\in SK_{l}$ if and only if there exists 
$X, Y, Z, T\in{\mathcal O}_r$ such that 
$$h_{\varpi c,d}wB_{i,\alpha,\beta}w h_{\varpi c,d}^{-1}=\begin{pmatrix}X&Y\\-\hat{\Delta}Y+\varpi^l Z&X+\varpi^l T\end{pmatrix}.$$
In this hypothesis we have $2c\varpi^{i+1}=\varpi^l T,$ i.e $c=\varpi^{l-i-1}e, e \in{\mathcal O}_r.$ We now look for a necessary condition on $d.$
We also have 
\begin{equation}-\hat{\Delta} Y\varpi^i+\varpi^l Z=(\varpi^{i+j+1}\beta'-\varpi^{r-i}e^2)/d.\label{constraint}
\end{equation}
From this last equation we have to distinguish two cases: $i+j+1<l$ or $i+j+1\geq l.$

$\Diamond$$\Diamond$In the first case $i+j+1<l.$ 
Because $Y=d\varpi^i$,
$-\hat{\Delta} Y$ is necessarily of valuation $i+k$, but from the equation (\ref{constraint}) and the inequality $i+j+1<l,$  we must have $k=j+1.$
We therefore have that if $k\not=j+1$ then $h_{\varpi c,d}wB_{i,\alpha,\beta}w h_{\varpi c,d}^{-1}$ is not in $SK_l$. The character of the representation is zero.

If $k=j+1$ then we have to solve the equation 
\begin{equation}\label{equationd}
-\hat{\Delta}' d\varpi^{i+k}+\varpi^l Z=(\varpi^{i+j+1}\beta'-\varpi^{r-i}e^2)/d.
\end{equation}
There are two cases to consider. 
This equation modulo $\varpi^l$ gives:
$$-\hat{\Delta}' d\varpi^{i+k}=\varpi^{i+j+1}\beta'/d\; \text{mod} \;\varpi^l.$$

The first case is when $-\frac{\hat{\Delta}'}{\beta'}$ is not a square, then   there is no solution in $d$  to this equation implying  $h_{\varpi c,d}wB_{i,\alpha,\beta}w h_{\varpi c,d}^{-1}$ is not in $SK_l$. 
The value of the character is therefore equal to zero.

The second case is when $-\frac{\hat{\Delta}'}{\beta'}=\Gamma^{-2}$, the equation (\ref{equationd}) implies that $d=\pm \Gamma+\varpi^{l-i-k}f$ with $f\in {\mathcal O}_r.$
As a result we obtain

\begin{equation*}
h_{\varpi c,d}wB_{i,\alpha,\beta}w h_{\varpi c,d}^{-1}=\begin{pmatrix}\alpha-\varpi^{l}e& d\varpi^i\\
(\varpi^{i+k}\beta'-\varpi^{2l-i}e^2)/d& \alpha+\varpi^{l}e
\end{pmatrix}=sk\in SK_l
\end{equation*}
with 
\begin{eqnarray*}
&&s =\begin{pmatrix}\alpha&d\varpi^i\\-d\varpi^i\hat{\Delta}& \alpha
\end{pmatrix} \\
&&k=\frac{1}{\alpha^2+d^2\varpi^{2i}\hat{\Delta}}\begin{pmatrix}\alpha^2-\alpha\varpi^le-\varpi^{2i+k}\beta'&-d\varpi^{l+i}e
\\\alpha d\varpi^{i+k}\hat{\Delta}'-d\varpi^{i+k+l}\hat{\Delta}'e+\alpha d^{-1}(\varpi^{i+k}\beta'-\varpi^{2l-i}e^2)& \alpha^2+\alpha\varpi^l e+
\varpi^{2i+k}d^2\hat{\Delta}'
\end{pmatrix}.
\end{eqnarray*}
(one can very that indeed $k$ belongs to $K_l$.)
Therefore the character of  ${\Xi}_{\Delta,\theta} $ evaluated on the conjugacy class is equal to 
\begin{eqnarray*}
\chi_3=\sum_{e\in{\mathbb O}_{i} \atop d=\pm\Gamma+\varpi^{l-i-k}f, f\in{\mathbb O}_{i+k}}\theta(\begin{pmatrix}\alpha&d\varpi^i\\-d\varpi^i\hat{\Delta}& \alpha
\end{pmatrix})\psi((d\varpi^{i+k}\hat{\Delta}'+d^{-1}(\varpi^{i+k}\beta'-\varpi^{2l-i}e^2))\frac{\alpha}{\alpha^2+d^2\varpi^{2i}\hat{\Delta}}).
\end{eqnarray*}
We have not been able to simplify this formula further.

$\Diamond$$\Diamond$In the second case $i+j+1\geq l.$  In order for $h_{\varpi c,d}wB_{i,\alpha,\beta}w h_{\varpi c,d}^{-1}$ to belong to $SK_{l}$, we necessarily have
 $c=\varpi^{l-i-1}e, $ and there must exists $Z\in {\mathcal O}_r$  such that $-\hat{\Delta}' d\varpi^{i+k}+\varpi^l z=(\varpi^{i+j+1}\beta'-\varpi^{r-i}e^2)/d.$
We have to distinguish 2 cases:

 $\Diamond\Diamond\Diamond$ $i+k< l.$ There is no solution $Z$ to the previous equation because $v(\hat{\Delta}' d\varpi^{i+k}+\varpi^{i+j+1}\beta'-\varpi^{r-i}e^2)/d)=i+k.$  Therefore the value of the character is  zero on the conjugacy class.

$\Diamond\Diamond\Diamond$ $i+k\geq l$ Then $h_{\varpi c,d}wB_{i,\alpha,\beta}w h_{\varpi c,d}^{-1}$  belong to $SK_{l}$ if $c=\varpi^{l-i-1}e.$
We proceed analogously as in the previous case

\begin{equation*}
 h_{\varpi c,d}wB_{i,\alpha,\beta}w h_{\varpi c,d}^{-1}=\begin{pmatrix}\alpha-\varpi^{l}e& d\varpi^i\\
(\varpi^{i+j+1}\beta'-\varpi^{2l-i}e^2)/d& \alpha+\varpi^{l}e
\end{pmatrix}=sk\in SK_l
\end{equation*}with 
\begin{eqnarray*}
&& s =\begin{pmatrix}\alpha&d\varpi^i\\-d\varpi^i\hat{\Delta}& \alpha
\end{pmatrix}\\
&&k=\frac{1}{\alpha^2+d^2\varpi^{2i}\hat{\Delta}}\begin{pmatrix}\alpha^2-\alpha\varpi^le-\varpi^{2i+j+1}\beta'&-d\varpi^{l+i}e
\\\alpha(d\varpi^{i+k}\hat{\Delta}'+d^{-1}(\varpi^{i+j+1}\beta'-\varpi^{2l-i}e^2))& \alpha^2+\alpha\varpi^l e+
\varpi^{2i+k}d^2\hat{\Delta}'
\end{pmatrix}.
\end{eqnarray*}

Therefore the character of  ${\Xi}_{\Delta,\theta} $ evaluated on the conjugacy class is equal to 
\begin{eqnarray*}
\sum_{e\in{\mathbb O}_i\atop d\in{\mathbb O}_l^\times}\theta(\begin{pmatrix}\alpha&d\varpi^i\\-d\varpi^i\hat{\Delta}& \alpha
\end{pmatrix})\psi((d\varpi^{i+k}\hat{\Delta}'+d^{-1}(\varpi^{i+j+1}\beta'-\varpi^{2l-i}e^2))\frac{\alpha}{\alpha^2+d^2\varpi^{2i}\hat{\Delta}}).
\end{eqnarray*}
Noting  that $\frac{\alpha}{\alpha^2+d^2\varpi^{2i}\hat{\Delta}}=\frac{1}{\alpha}(1-\frac{\hat{\Delta}'}{\alpha^2}\varpi^{2i+k}d^2)$,
the value of the character is:
\begin{eqnarray*}
\chi_4=\sum_{e\in{\mathbb O}_i\atop d\in{\mathbb O}_l^\times}\theta(\begin{pmatrix}\alpha&d\varpi^i\\-d\varpi^i\hat{\Delta}& \alpha
\end{pmatrix})\psi(\frac{1}{\alpha}(d\varpi^{i+k}\hat{\Delta}'+d^{-1}(\varpi^{i+j+1}\beta'-\varpi^{2l-i}e^2))).
\end{eqnarray*}

We have not been able to simplify this formula further.

To summarize:
\medskip

\begin{tabular}{|c|c|c|c|c|}\hline
  &$I_{\alpha}$   & $D_{i,\alpha,\delta}$& $C_{i,\alpha,\beta}$ &$B_{i,\alpha,\beta}$ \\ 
\hline\hline
$tr({{\Xi}_{a,\Delta,\theta}}) (.)$ & $\begin{array}{c}(q^2-1)q^{r-2}\tilde{\chi}_a(\alpha^2)\\ \sigma(\alpha)\end{array}$
 &$\begin{array}{c}\delta_{i,r-1}\tilde{\chi}_a(\alpha)\tilde{\chi}_a(\delta)\\\sigma(\delta)(q-1)q^{r-2}\end{array}$ &$\begin{array}{c}-\delta_{i,r-1}\tilde{\chi}_a(\alpha^2)\\\sigma(\alpha)q^{r-2}(q+1)\end{array}$&${\text{Many  cases }\atop\text{No "simple" formula}}$ \\
\hline
\end{tabular}
\bigskip




\begin{appendix}
\section{Gauss sums, Kloosterman sums, Sali\'e sums}
We use the notations of section $3$.

Let $\lambda:({\mathcal O}_k,+)\rightarrow {\mathbb C}^\times$ a primitive  character, let $a\in  {\mathcal O}_k$,  we will denote  the 
quadratic Gauss sum $G(a, \lambda)$ to be:
\begin{equation*}
G(a,\lambda)=\sum_{x\in {\mathcal O}_k}\lambda(ax^2).
\end{equation*}

Let $a,b \in  {\mathcal O}_k$, the Kloosterman sum $K(a,b,\lambda)$ is defined as:
\begin{equation*}
K(a,b,\lambda)=\sum_{x\in {\mathcal O}_k^{\times }}\lambda(ax+bx^{-1}).
\end{equation*}

Let $\rho$ be a multiplicative character ${\mathcal O}_k^\times \rightarrow {\mathbb C}^\times$, the twisted Kloosterman sum 
$K(a,b,\lambda,\rho)$ is defined as:
\begin{equation*}
K(a,b,\lambda,\rho)=\sum_{x\in {\mathcal O}_k^{\times }}\rho(x)\lambda(ax+bx^{-1}).
\end{equation*}

Important example of twisted Kloosterman sum which appear in our work is the Sali\'e sum
 $S(a,b,\lambda)$  defined as:
\begin{equation*}
S(a,b,\lambda)=\sum_{x\in {\mathcal O}_k^{\times }}(\frac{x}{{\mathcal O}_k})\lambda(ax+bx^{-1}),
\end{equation*}
 where $ (\frac{x}{{\mathcal O}_k})$ denotes the Legendre symbol in ${\mathcal O}_k$,  which is defined for every $x\in {\mathcal O}_k$ and is equal to
\begin{equation*} (\frac{x}{{\mathcal O}_k})=
\left\{
\begin{array}{ll}0 &\text{if}\; x \;\text{is not invertible}\\1  &\text{if}\; x \;\text{is a square} \\-1  \;&\text{otherwise}.
\end{array}
\right.
\end{equation*}
Note that the Legendre symbol restricted to ${\mathcal O}_k^{\times}$ is a group morphism with value in $\{+1,-1\}$ which factor through the group   
$\Bbbk^\times.$ 

Remark: In order to keep track of the dependence of $k,$ we will sometimes use the notation $G_k, K_k, S_k$ for the Gauss, Kloosterman, Sali\'e sum associated to ${\mathcal O}_k.$

In \cite{Sz} quadratic Gauss sum are studied and computed for any finite commutative ring of odd characteristic. We apply his results to the case  of the ring 
${\mathcal O}_k.$ With his notations, we have $d_{{\mathcal O}_k}=k$, and the theorem 6.2 of \cite{Sz} can be stated as:
\begin{proposition}
\begin{eqnarray*}
&&G(1,\lambda)^2= (\frac{-1}{{\mathcal O}_k})^k q^k\\
&&G(ab,\lambda)= (\frac{a}{{\mathcal O}_k})^k G(b,\lambda), a,b\in {\mathcal O}_k^{\times}.
\end{eqnarray*}
\end{proposition}

From this theorem we obtain the following result which is needed for the evaluation of the characters of cuspidal representations for conjugacy class of type C.

\begin{proposition}\label{KloostermanGauss}
Let $\lambda:({\mathcal O}_k,+)\rightarrow {\mathbb C}^\times $ be a primitive character, let $\eta\in {\mathcal O}_k$ an invertible element which is not a square, we have:
\begin{equation*}
\sum_{e,f\in{ {\mathcal O}_k}} \lambda(e^2-\eta f^2 )=(-q)^k.
\end{equation*}
\end{proposition}
\proof 
Let $S=\sum_{e,f\in{ {\mathcal O}_k}} \lambda(e^2-\eta f^2 )=G(1,\lambda)G(-\eta,\lambda).$
We have $S=G(1,\lambda)G(-\eta,\lambda)=G(1,\lambda)(\frac{-\eta}{{\mathcal O}_k})^k G(1,\lambda)=
(\frac{-\eta}{{\mathcal O}_k})^k (\frac{-1}{{\mathcal O}_k})^k q^k=(\frac{\eta}{{\mathcal O}_k})^k  q^k=(-q)^k.$
\endproof

In the evaluation of characters of cuspidal representations, one needs an explicit expression for $T(b,\eta,\lambda)=\sum_{c\in {\mathcal O}_k, d\in {\mathcal O}_k^{\times}}
\lambda(d^{-1}(b-c^2)+d \eta)$ where $b \in {\mathcal O}_k, \eta\in   {\mathcal O}_k^\times$, $\eta$ not a square and $\lambda:({\mathcal O}_k,+)\rightarrow {\mathbb C}^\times $ is  a primitive character. This sum is a twisted Kloosterman sum.
Indeed we have:
\begin{eqnarray*}
&&T(b,\eta,\lambda)=\sum_{c\in {\mathcal O}_k, d\in {\mathcal O}_k^{\times}}
\lambda(d^{-1}(b-c^2)+d \eta)\\
&&=\sum_{d\in {\mathcal O}_k^{\times}}
\lambda(db+d^{-1}\eta)G(-d,\lambda)\\
&&=\sum_{d\in {\mathcal O}_k^{\times}}
\lambda(db\eta+d^{-1} )G(-d\eta,\lambda)\\
&&=G(-1,\lambda)\sum_{d\in {\mathcal O}_k^{\times}}(\frac{d\eta}{{\mathcal O}_k})^k
\lambda(db\eta+d^{-1} )\\
&&=G(-1,\lambda)(-1)^k\sum_{d\in {\mathcal O}_k^{\times}}(\frac{d}{{\mathcal O}_k})^k
\lambda(db\eta+d^{-1} ).
\end{eqnarray*}
As a result when $k$ is even we get a Kloosterman sum and when $k$ is odd we obtain a Sali\'e sum.
The following result give a simple formula of the evaluation of this sum for any $k$. When $k$ is even, and ${\mathcal O}_k={\mathbb Z}/p^k{\mathbb Z},$ this is the classical formula for evaluation of Kloosterman sum obtained by  H.Sali\'e in 1931. 
In the case where ${\mathcal O}$ is the ring of integer of a $p$-adic field $F$, we could obtain the evaluation of  these sums by applying the results of \cite{Pa1,Pa2} to a number field having $F$ at some place.
This is not completely direct and do not cover the case where the local field is of positive characteristic,  we prefer to give a direct proof of it using a generalization of the method of \cite{LCW}.

\begin{proposition}

Let $\lambda:({\mathcal O}_k,+)\rightarrow {\mathbb C}^\times $ be a primitive character, let $\eta\in {\mathcal O}_k$ an invertible element which is not a square,  we have:

\begin{equation*}\label{evaluationkloosterman1}
\sum_{c\in {\mathcal O}_k, d\in {\mathcal O}_k^{\times}}
\lambda(d^{-1}(b-c^2)+d \eta)=\left\{
\begin{array}{ll}(-q)^k(\lambda(2u)+\lambda(-2u)), &\text{if}  \;u^2=\eta b\;\; \text{is invertible},\\-q\; &\text{if}\; k=1 \;\text{and}\; b=0,\\0 \;&\text{otherwise}.
\end{array}
\right.
\end{equation*}

\end{proposition}
\proof

The sum $T(b,\eta,\lambda)$ can be expressed as:

\begin{eqnarray*}
&&T(b,\eta,\lambda)=\sum_{c\in {\mathcal O}_k, d\in {\mathcal O}_k^{\times}}
\lambda(d^{-1}(b-c^2)+d \eta)\\
&&=\sum_{x\in {\mathcal O}_k}
\lambda(x)\vert\{(c,d)\in {\mathcal O}_k  \times  {\mathcal O}_k^{\times}, x=\eta d+(b-c^2) d^{-1}\}\vert \\
&&=\sum_{x\in {\mathcal O}_k}
\lambda(x)\vert E(x)\cap  ({\mathcal O}_k  \times  {\mathcal O}_k^{\times})\vert 
\end{eqnarray*}
where $E(x)=\{(c,d)\in {\mathcal O}_k  \times  {\mathcal O}_k, x d=\eta d^2+(b-c^2) \}.$

Noting that $E(x)\cap  ({\mathcal O}_k  \times  \varpi {\mathcal O}_k)=
E(x+\varpi^{k-1})\cap  ({\mathcal O}_k  \times  \varpi {\mathcal O}_k),$ we obtain that 
\begin{eqnarray*}
&&\sum_{x\in {\mathcal O}_k}
\lambda(x)\vert E(x)\cap  ({\mathcal O}_k  \times  \varpi{\mathcal O}_k)\vert =
\sum_{x\in {\mathcal O}_k}
\lambda(x)\vert E(x+\varpi^{k-1})\cap  ({\mathcal O}_k  \times 
 \varpi{\mathcal O}_k)\vert\\
&&=\sum_{x\in {\mathcal O}_k}
\lambda(x-\varpi^{k-1})\vert E(x)\cap  ({\mathcal O}_k  \times 
 \varpi{\mathcal O}_k)\vert\\
&&=\lambda(-\varpi^{k-1})\sum_{x\in {\mathcal O}_k}
\lambda(x)\vert E(x)\cap  ({\mathcal O}_k  \times  \varpi{\mathcal O}_k)\vert\\
&&=0.\end{eqnarray*}
As a result we obtain:
$T(b,\eta,\lambda)=\sum_{x\in {\mathcal O}_k}
\lambda(x)\vert E(x)\vert.$

Let $F(x)=\{(c,d)\in {\mathcal O}_k^{ 2}, d^2-\eta c^2=x\},$ a simple computation shows that $E(x)=F((\frac{x}{2})^2-\eta b).$
As a result, if we introduce $\rho:{\mathcal O}_k\rightarrow {\mathbb N}, \rho(y)=
\vert F(y) \vert,$ and noting that $\rho(y)=\rho(u^2 y)$ if $u$ is invertible, 
 $$T(b,\eta,\lambda)=\sum_{x\in {\mathcal O}_k}\lambda(x)
\rho(x^2-4\eta b).$$

Using a straighforward generalization of the argument of \cite{LCW}, $\rho$ can be evaluated exactly and is a function of the valuation
$$
\rho(y)
=\left\{
\begin{array}{ll}q^{2(k-\lfloor \frac{k+1}{2} \rfloor)} &\text{if} \; y=0\\
(q+1)q^{k-1-v(y)}&\text{if} \;v(y)\; \text{even (including 0)}\\
0 \;&\text{if} \;v(y)\;\text{odd}.
\end{array}
\right.
$$

Let us recall the argument of \cite{LCW} generalized in our setting.
We consider the different cases.
\begin{enumerate}
\item $b$ invertible and is a square
\item $b$ invertible and is a non-square
\item $b$ non invertible
\end{enumerate}
In the case $1)$, $\eta b$ is not a square, this is also equivalent by Hensel lemma to the fact that it is not a square in the residual field. Therefore $x^2-4\eta b$ is invertible for all $x\in {\mathcal O}_k$, because otherwise it would vanish in the residual field contradicting that $\eta b$ is not a square.
Therefore
\begin{eqnarray*}
&&T(b,\eta,\lambda)=\sum_{x\in {\mathcal O}_k}\lambda(x)
\rho(x^2-4\eta b)\\
&&=\sum_{x\in {\mathcal O}_k}\lambda(x)(q+1)q^{k-1}=0.
\end{eqnarray*}
In the case $3)$ $b$ belongs to $\varpi{\mathcal O}_k.$
When $k\geq 2$, we have

 \begin{eqnarray*}
&&T(b,\eta,\lambda)=\sum_{x\in {\mathcal O}_k^\times }\lambda(x)
\rho(x^2-4\eta b)+\sum_{x\in\varpi {\mathcal O}_k }\lambda(x)
\rho(x^2-4\eta b)\\
&&=\sum_{x\in {\mathcal O}_k^\times }\lambda(x)
(q+1)q^{k-1}+\sum_{x\in\varpi {\mathcal O}_k }\lambda(x)
\rho(x^2-4\eta b)
\end{eqnarray*}
The first sum is $0$ after having used 
 the following property, direct generalisation of the lemma 5.1 of \cite{LCW}
\begin{equation}\label{sumlambda}
\sum_{x\in \varpi^j {\mathcal O}_k^\times } \lambda(x)
=\left\{
\begin{array}{ll}0 &\text{if} \; j<k-1\\
-1&\text{if} \;j=k-1\\
1&\text{if} \;j=k.
\end{array}
\right.
\end{equation}
To evaluate the second sum we notice that $(x+\varpi^{k-1})^2=x^2$ when $x\in\varpi {\mathcal O}_k$, hence:
\begin{eqnarray*}
&&\sum_{x\in\varpi {\mathcal O}_k }\lambda(x)
\rho(x^2-4\eta b)=\sum_{x\in\varpi {\mathcal O}_k }\lambda(x)
\rho((x+\varpi^{k-1})^2-4\eta b)\\
&&=\sum_{x\in\varpi {\mathcal O}_k }\lambda(x-\varpi^{k-1})
\rho(x^2-4\eta b)=\lambda(-\varpi^{k-1})\sum_{x\in\varpi {\mathcal O}_k }\lambda(x)
\rho(x^2-4\eta b),
\end{eqnarray*}
implying the vanishing of the second sum. Therefore $T(b,\eta,\lambda)=0.$

Note that when $k=1$, we necessarily have $b=0,$   and 
\begin{eqnarray*}
&&T(b,\eta,\lambda)=\sum_{x\in {\mathcal O}_1}\lambda(x)
\rho(x^2)\\
&&=\lambda(0)\rho(0)+\sum_{x\in \Bbbk^\times }\lambda(x)(q+1)\\
&&=1-(q+1)=-q.
\end{eqnarray*}

In the remaining case $2)$, we have $\epsilon b=u^2$ with $u$ invertible, therefore 
$  T(b,\eta,\lambda)=\sum_{x\in {\mathcal O}_k }\lambda(x)
\rho((x-2u)(x+2u)).$

$(x-2u)(x+2u)$ is non invertible if and only if $x-2u$ or $x+2u$ has a strictly positive valuation.
We denote $X_j^\pm=\pm 2u+\varpi^j {\mathcal O}_k, j\geq 1.$
We have 
\begin{eqnarray*}
&&T(b,\eta,\lambda)=\sum_{x\in {\mathcal O}_k }\lambda(x)
\rho((x-2u)(x+2u))\\
&&=\sum_{x\in {\mathcal O}_k \setminus{\cup_{j=1}^k X_j^\pm}}\lambda(x)\rho(1)+
\sum_{\epsilon=\pm}\sum_{j=1}^k\sum_{x \in X_j^\epsilon}\lambda(x)
\rho(\varpi^j)\\
&&=\sum_{x\in {\mathcal O}_k}\lambda(x)\rho(1)+
\sum_{\epsilon=\pm}\sum_{j=1}^k\sum_{x \in X_j^\epsilon}\lambda(x)
(\rho(\varpi^j)-\rho(1))\\
&&=\sum_{x\in {\mathcal O}_k}\lambda(x)\rho(1)+
\sum_{\epsilon=\pm}\sum_{j=1}^k\sum_{x \in \varpi^j{\mathcal O}_k }\lambda(\epsilon u)\lambda(x)
(\rho(\varpi^j)-\rho(1)).
\end{eqnarray*}
Using $\sum_{x\in \varpi^j{\mathcal O}_k}\lambda(x)=0$  for $j=0,...,k-2$ we obtain
\begin{eqnarray*}
&&T(b,\eta,\lambda)=
\sum_{j=k-1}^k\sum_{x \in\varpi^j{\mathcal O}_k  }(\lambda(2u)+\lambda(-2u))\lambda(x)
(\rho(\varpi^j)-\rho(1)).
\end{eqnarray*}
Applying (\ref{sumlambda}) we end up with:
\begin{eqnarray*}
&&T(b,\eta,\lambda)=
-(\lambda(2u)+\lambda(-2u))(\rho(\varpi^{k-1})-\rho(1))+
(\lambda(2u)+\lambda(-2u))(\rho(0)-\rho(1))\\
&&=(\lambda(2u)+\lambda(-2u))(\rho(0)-\rho(\varpi^{k-1}))\\
&&=(-q)^k ((\lambda(2u)+\lambda(-2u))).\end{eqnarray*}
This ends the proof of the proposition.

\endproof

The rest of this section is devoted to the evaluation of the sum $$\chi_2(i,j,k)=q^{2(i-l)}\theta(I_{\alpha})\sum_{ c\in {\mathbb O}_{r-i-1}, d\in {\mathbb O}_{r-i}^{\times}}
\psi(\frac{\varpi^i}{\alpha}(\frac{\varpi^{j+1}\beta'-\varpi^{2}c^2}{d}-d\varpi^{j+k}\hat{\Delta}'\beta'))$$ with the conditions 
$l\leq i<r-1$, $0\leq j\leq r-1$, $1\leq k\leq r, 1+i+j\leq r.$
This sum is the character of ${\Xi}_{\Delta,\theta}$ evaluated on the conjugacy class $B_{i,\alpha,\beta}$ when $l\leq i<r-1.$

Precise evaluation of these kind of sums have been given by Maeda in \cite{Ma} for the case ${\mathcal O}={\mathbb Z}_p.$ We will show that $\chi_2(i,j,k)$
 can always be expressed in term of  Gauss Sums, Kloosterman sums and Sali\'e  sums. In most cases one can further evaluate them but in the case where $k=1$ there are cases where the evaluation amount to evaluate Kloosterman sums in the case where there is no closed formula for them.

\begin{eqnarray*}
\chi_2(i,j,k)&=&q^{2(i-l)-1}\theta(I_{\alpha})\sum_{ c\in {\mathbb O}_{r-i}, d\in {\mathbb O}_{r-i}^{\times}}
\psi(\frac{\varpi^i}{\alpha}(\frac{\varpi^{j+1}\beta'-\varpi^{2}c^2}{d}-d\varpi^{j+k}\hat{\Delta}'\beta'))\\
&=&q^{2(i-l)-1}\theta(I_{\alpha})\sum_{ c\in {\mathcal O}_{r-i}, d\in {\mathcal O}_{r-i}^{\times}}
\lambda(\frac{\varpi^{j+1}\beta'-\varpi^{2}c^2}{d}-d\varpi^{j+k}\hat{\Delta}'\beta')\\
&=&q^{2(i-l)-1}\theta(I_{\alpha})\sum_{ c\in {\mathcal O}_{r-i}, d\in {\mathcal O}_{r-i}^{\times}}
\lambda((\varpi^{j+1}\beta'-\varpi^{2}c^2)d-\frac{\varpi^{j+k}\hat{\Delta}'\beta'}{d})
\end{eqnarray*}
where $\lambda:{\mathcal O}_{r-i}\rightarrow {\mathbb C}^{\times}$ is the primitive character factor map of the character  $z\mapsto\psi(\varpi^i\frac{z}{\alpha}).$

When $i=r-2$ then the elements on which $\lambda$ is evaluated are all $0.$ Therefore we obtain 
$\chi_2(i,j,k)=q^{2(i-l)-1}\vert{\mathcal O}_{2}\vert
\vert{\mathcal O}_{2}^\times\vert=q^{r-2}(q-1).$
We now assume $i<r-2.$

We have to distinguish two cases: $j\geq1$ and $j=0$

$\Diamond$
If $j\geq 1$

Denoting $\mu:{\mathcal O}_{r-i-2}\rightarrow {\mathbb C}^{\times}$ the primitive character factor map of the character  $z\to\lambda(\varpi^2 z),$ 
we obtain:
\begin{eqnarray*}
\chi_2(i,j,k)&=&q^{2(i-l)+3}\theta(I_{\alpha})\sum_{ c\in {\mathcal O}_{r-i-2}\atop d\in {\mathcal O}_{r-i-2}^{\times}}
\mu((\varpi^{j-1}\beta'-c^2)d-\frac{\varpi^{j-1+k-1}\hat{\Delta}'\beta'}{d})\\
&=&q^{2(i-l)+3}\theta(I_{\alpha})\sum_{d\in {\mathcal O}_{r-i-2}^{\times}}
\mu(\varpi^{j-1}\beta'd-\frac{\varpi^{j+k-2}\hat{\Delta}'\beta'}{d})G_{r-i-2}(\mu,-d)\\
&=&q^{2(i-l)+3}\theta(I_{\alpha})G_{r-i-2}(\mu,-1)\sum_{d\in {\mathcal O}_{r-i-2}^{\times}}
\mu(\varpi^{j-1}\beta'd-\frac{\varpi^{j+k-2}\hat{\Delta}'\beta'}{d})(\frac{d}{{\mathcal O}_{r-i-2}})^{r-i-2}.
\end{eqnarray*}
When $i+j+1=r$, we have 
\begin{eqnarray*}
\chi_2(i,j,k)&&=q^{2(i-l)+3}\theta(I_{\alpha})G_{r-i-2}(\mu,-1)\sum_{d\in {\mathcal O}_{r-i-2}^{\times}}
(\frac{d}{{\mathcal O}_{r-i-2}})^{r-i-2}\\
&=&q^{2(i-l)+3}\theta(I_{\alpha})G_{r-i-2}(\mu,-1)\times \left\{
\begin{array}{ll}\vert {\mathcal O}_{r-i-2}^{\times}\vert &\text{if} \;i \;\text{is even}\\
0  &\text{if}\; i \;\text{is odd} .
\end{array}
\right.
\end{eqnarray*}

When $i+j+1<r$ we 
 denote $\tilde{\mu}:{\mathcal O}_{r-i-j-1}\rightarrow {\mathbb C}^{\times}$ the primitive character factor map of the character  $z\to\mu(\varpi^{j-1} z),$
\begin{eqnarray*}
\chi_2(i,j,k)&=&q^{2(i-l)+3}\theta(I_{\alpha})G_{r-i-2}(\mu,-1)\times \left\{
\begin{array}{ll}K_{r-i-j-1}(\beta',-\varpi^{k-1}\hat{\Delta}'\beta',\tilde{\mu}) &\text{if} \;i \;\text{is even}\\
S_{r-i-j-1}(\beta',-\varpi^{k-1}\hat{\Delta}'\beta',\tilde{\mu})  &\text{if}\; i \;\text{is odd} .
\end{array}
\right.
\end{eqnarray*}

When $k>1$ we can further simplify these expressions. Indeed using the fact that the map $d\mapsto d -\frac{\varpi^{k-1}\hat{\Delta}'}{d}$ is a bijection from 
${\mathcal O}_{r-i-j-1}^\times $to $ {\mathcal O}_{r-i-j-1}^\times$ when $k>1$, we can write 

\begin{eqnarray*}
\chi_2(i,j,k)&=&q^{2(i-l)+3}\theta(I_{\alpha})G_{r-i-2}(\mu,-1)\times \left\{
\begin{array}{ll}K_{r-i-j-1}(\beta',0,\tilde{\mu}) &\text{if} \;i \;\text{is even}\\
S_{r-i-j-1}(\beta',0,\tilde{\mu})  &\text{if}\; i \;\text{is odd} .
\end{array}
\right.
\end{eqnarray*}

When $k=1$; we have
\begin{eqnarray*}
\chi_2(i,j,k)&=&q^{2(i-l)+3}\theta(I_{\alpha})G_{r-i-2}(\mu,-1)\times \left\{
\begin{array}{ll}K_{r-i-j-1}(\beta',-\hat{\Delta}'\beta',\tilde{\mu}) &\text{if} \;i \;\text{is even}\\
S_{r-i-j-1}(\beta',-\hat{\Delta}'\beta',\tilde{\mu})  &\text{if}\; i \;\text{is odd} .
\end{array}
\right.
\end{eqnarray*}
Note that  only  the case $r-i-j-1=1$ and $i$ even cannot be further simplified.

$\Diamond$ If $j=0$ we have

\begin{eqnarray*}
\chi_2(i,0,k)&=&q^{2(i-l)-1}\theta(I_{\alpha})\sum_{ c\in {\mathcal O}_{r-i}\atop d\in {\mathcal O}_{r-i}^{\times}}
\lambda((\varpi\beta'-\varpi^{2}c^2)d-\frac{\varpi^{k}\hat{\Delta}'\beta'}{d})\\
&=&q^{2(i-l)-1}q^2\theta(I_{\alpha})\sum_{ c\in {\mathcal O}_{r-i-1}\atop d\in {\mathcal O}_{r-i-1}^{\times}}
\tilde{\lambda}((\beta'-\varpi^{}c^2)d-\frac{\varpi^{k-1}\hat{\Delta}'\beta'}{d})
\end{eqnarray*}
with  $\tilde{\lambda}: {\mathcal O}_{r-i-1}\to {\mathbb C}$ primitive character factor map of  $z\mapsto\lambda(\varpi z).$
If we still denote $ \mu:{\mathcal O}_{r-i-2}\rightarrow {\mathbb C}^{\times}$ the primitive character factor map of  the character  $z\to\tilde{\lambda}(\varpi z),$ 
we obtain 
\begin{eqnarray*}
\chi_2(i,0,k)&=&q^{2(i-l)-1}q^2q\theta(I_{\alpha})\sum_{ c\in {\mathcal O}_{r-i-2}\atop d\in {\mathcal O}_{r-i-1}^{\times}}
\tilde{\lambda}(\beta'd-\frac{\varpi^{k-1}\hat{\Delta}'\beta'}{d})\mu(-p(d)c^2)\\
&=&
q^{2(i-l)+2}\theta(I_{\alpha})\sum_{ d\in {\mathcal O}_{r-i-1}^{\times}}
\tilde{\lambda}(\beta'd-\frac{\varpi^{k-1}\hat{\Delta}'\beta'}{d})G_{r-i-2}(\mu,-p(d)).
\end{eqnarray*}
where $p(d)$ is the projection of $d$ in  ${\mathcal O}_{r-i-2}.$
As a result we get:
\begin{eqnarray*}
\chi_2(i,0,k)&=&q^{2(i-l)+2}\theta(I_{\alpha})G_{r-i-2}(\mu,-1) \sum_{ d\in {\mathcal O}_{r-i-1}^{\times}}
\tilde{\lambda}(\beta'd-\frac{\varpi^{k-1}\hat{\Delta}'\beta'}{d})(\frac{p(d)}{{\mathcal O}_{r-i-2}})^{r-i-2}\\
&=&q^{2(i-l)+2}\theta(I_{\alpha})G_{r-i-2}(\mu,-1) \sum_{ d\in {\mathcal O}_{r-i-1}^{\times}}
\tilde{\lambda}(\beta'd-\frac{\varpi^{k-1}\hat{\Delta}'\beta'}{d})(\frac{d}{{\mathcal O}_{r-i-1}})^{r-i-2}\\
&=&q^{2(i-l)+2}\theta(I_{\alpha})G_{r-i-2}(\mu,-1) \times \left\{
\begin{array}{ll}K_{r-i-1}(\beta',-\varpi^{k-1}\hat{\Delta}'\beta',\tilde{\lambda}) &\text{if} \;i \;\text{is even}\\
S_{r-i-1}(\beta',-\varpi^{k-1}\hat{\Delta}'\beta',\tilde{\lambda})  &\text{if}\; i \;\text{is odd} .
\end{array}
\right.
\end{eqnarray*}

When $k>1$ we can further simplify these expressions. Indeed using the fact that the map $d\mapsto d -\frac{\varpi^{k-1}\hat{\Delta}'}{d}$ is a bijection from 
${\mathcal O}_{r-i-1}^\times $to $ {\mathcal O}_{r-i-1}^\times$ when $k>1$, we can write 

\begin{eqnarray*}
\chi_2(i,0,k)&=&q^{2(i-l)+2}\theta(I_{\alpha})G_{r-i-2}(\mu,-1)\times \left\{
\begin{array}{ll}K_{r-i-1}(\beta',0,\tilde{\lambda}) &\text{if} \;i \;\text{is even}\\
S_{r-i-1}(\beta',0,\tilde{\lambda})  &\text{if}\; i \;\text{is odd} .
\end{array}
\right.
\end{eqnarray*}

When $k=1$; we have
\begin{eqnarray*}
\chi_2(i,0,1)&=&q^{2(i-l)+2}\theta(I_{\alpha})G_{r-i-2}(\mu,-1)\times \left\{
\begin{array}{ll}K_{r-i-1}(\beta',-\hat{\Delta}'\beta',\tilde{\lambda}) &\text{if} \;i \;\text{is even}\\
S_{r-i-1}(\beta',-\hat{\Delta}'\beta',\tilde{\lambda})  &\text{if}\; i \;\text{is odd} .
\end{array}
\right.
\end{eqnarray*}
Note that  only  the case $i=r-2$  cannot be further simplified.

\end{appendix}
\medskip
{\bf Acknowledgments}
I am grateful  to I.Badulescu for numerous discussions.

\end{document}